\documentclass[number,citesort,seceqn,dvips]{arxbj}
\usepackage{upgreek}
\usepackage{graphicx}

\aid{0}
\volume{17}
\issue{2}
\pubyear{2011}
\firstpage{643}
\lastpage{670}
\doi{10.3150/10-BEJ286}

\makeatletter

    \newtheorem{thrm}{Theorem}[section]
    \newtheorem{prop}[thrm]{Proposition}
    
    \newtheorem{cllry}[thrm]{Corollary}
    \newtheorem{lmma}[thrm]{Lemma}
    \newremark{remk}[thrm]{Remark}
    \newproclaim{defn}[thrm]{Definition}
    
    \newproclaim{cond}{Assumption}

\def\ed{{\stackrel{\frak {D}}{=}}}
\def\ld{{\stackrel{\frak {D}}{\longrightarrow}}}

    \def\bbl{{\mathbb L}}

    \def\calB{{\mathcal B}}
    \def\calF{{\mathcal F}}

    \def\calL{{\mathcal L}}

    \def\calL{{\mathcal L}}
    \def\calN{{\mathcal N}}

    \def\calS{{\mathcal S}}
    
    \def\bbr{{\mathbb R}}
    
    \def\bbe{{\mathbb E }}
    \def\bbp{{ \mathbb P }}

    \def\iid{ {\stackrel{\mathrm{i.i.d.}}{\sim}} }
    
\makeatother

\begin{document}
\begin{frontmatter}

\title{Sieve-based confidence intervals and bands for L\'evy densities}
\runtitle{Confidence intervals and bands for L\'evy densities}

\begin{aug}
\author{\fnms{Jos\'e E.} \snm{Figueroa-L\'opez}\corref{}\ead[label=e1]{figueroa@purdue.edu}}

\runauthor{J. E. Figueroa-L\'opez}

\address{Department of Statistics,
Purdue University,
West Lafayette,
IN 47907,
USA.\\
\printead{e1}}
\end{aug}

\received{\smonth{11} \syear{2008}}
\revised{\smonth{5} \syear{2010}}

%
\begin{abstract}
The estimation of the L\'evy density,
the infinite-dimensional parameter controlling the
jump dynamics of a L\'evy process, is considered here
under a discrete-sampling scheme. In this setting,
the jumps are latent variables, the statistical
properties of which can be assessed when the frequency and time
horizon of
observations increase to infinity at suitable rates.
Nonparametric estimators for the L\'evy density
based on \textit{Grenander's method of sieves}
was proposed in Figueroa-L\'opez [\textit{IMS Lecture Notes} \textbf
{57} (2009) 117--146].
In this paper, central limit theorems for these sieve estimators, both
pointwise and uniform on
an interval away from the origin,
are obtained, leading to pointwise confidence intervals and bands for the
L\'evy density. In the pointwise case,
our estimators converge to the L\'evy density
at a rate that is arbitrarily close to the rate of the minimax risk of
estimation on smooth L\'evy densities.
In the case of uniform bands and discrete regular sampling, our
results are consistent with the case of density estimation, achieving a
rate of order arbitrarily close to $\log^{-1/2}(n)
\cdot n^{-1/3}$, where $n$ is the number of observations. The
convergence rates are valid,
provided that $s$ is smooth enough and that the time horizon
$T_{n}$ and the dimension of the sieve are appropriately chosen in
terms of $n$.

\end{abstract}

%

\end{frontmatter}
%
\section{Introduction}\label{sec1}

\subsection{Motivation and preliminary background}
In the past decade, L\'evy processes have received
a great deal of attention, fueled by numerous applications
in the area of mathematical finance,
to the extent that
L\'evy processes have become a fundamental building block
in the modeling of asset prices with jumps
(see, e.g., \cite{Cont2003} and \cite{Fig2010} for further
information about this field).
The simplest of these models postulates that
the price of a commodity (say a stock) at time $t$ is
given as an exponential function of a L\'evy process $X:=\{X_{t}\}
_{t\geq 0}$.
Even this simple extension of the classical
Black--Scholes model, in which $X$ is simply a Brownian motion with
drift, is able to account for
several fundamental empirical features commonly
observed in time series of asset returns,
such as heavy tails, high kurtosis
and asymmetry.
L\'evy processes, as models capturing some of
the most important features of returns and
as ``first-order approximations'' to other more
accurate models, are fundamental
for developing and testing successful
statistical methodologies. However, even in such
parsimonious models, there are several issues
concerning the performing of statistical inference by standard
likelihood-based methods.

A L\'evy process is the ``discontinuous sibling'' of a Brownian
motion. Concretely, $X=\{X_{t}\}_{t\geq 0}$
is a L\'evy process if $X$ has independent
and stationary increments, its paths are right-continuous
with left limits and it has no fixed jump times.
The later condition means that, for any $t>0$,
$
\bbp[\Delta X_{t}\neq 0]=0,
$
where $\Delta{X}_{t}:=X(t)-\lim_{s\nearrow t}X_{s}$
is the magnitude of the ``jump'' of $X$ at time $t$.
Any L\'evy process can be constructed from
the superposition of a Brownian motion with
drift, $\sigma W_{t}+bt$, a compound Poisson process and
the limit process resulting from making the jump
intensity of a compensated compound Poisson process,
$Y_{t}-\bbe Y_{t}$, go to
infinity while simultaneously allowing jumps
of smaller sizes.
Formally, $X$ admits a decomposition of the form
%
\begin{equation}\label{LevyItoDecmp}
X_{t}= bt+\sigma B_{t}+
\lim_{\varepsilon\searrow 0}
\int_{0}^{t}\!\!\int_{\varepsilon\leq|x|\leq1}
x
(\mu-\bar\mu)(\mathrm{d}x,\mathrm{d}s)+
\int_{0}^{t}\!\!\int_{|x|>1} x
\mu(\mathrm{d}x,\mathrm{d}s),
\end{equation}
where $B$ is a standard Brownian motion and
$\mu$ is an independent Poisson measure on
$\bbr_{+}\times\bbr\backslash\{0\}$ with
mean measure $\bar\mu(\mathrm{d}x,\mathrm{d}t):= \nu(\mathrm{d}x)\,\mathrm{d}t$.
Thus, L\'evy processes are determined by three parameters:
a nonnegative real $\sigma^{2}$, a real $b$ and
a measure $\nu$ on $\mathbb{R}\backslash\{0\}$
such that $\int(x^{2}\wedge1)\nu(\mathrm{d}x)<\infty$.
The measure $\nu$ controls the jump dynamics of the
process $X$, in that $\nu(A)$ gives the average number of jumps (per
unit time)
whose magnitudes fall in a given set $A\in\calB(\bbr)$.
A common assumption in L\'evy-based financial models is that
$\nu$ is determined by a function
$s\dvtx \bbr\backslash\{0\}\rightarrow[0,\infty)$, called
the \textit{L\'evy density}, as follows:
\[
\nu(A)=\int_{A} s(x) \,\mathrm{d}x   \qquad \forall
A\in\calB(\bbr\backslash\{0\}).
\]
Intuitively, the value of $s$ at $x_{0}$ provides
information on the frequency of jumps with sizes
``close'' to $x_{0}$.

\subsection{The statistical problem and methodology}
We are interested in estimating, in a nonparametric
fashion, the L\'evy density $s$ over a window of estimation
$D:=[a,b]\subset\bbr\backslash\{0\}$,
based on discrete observations of the process on a finite
interval $[0,T]$. In general, $s$ can blow up around the origin and, hence,
we consider only domains~$D$ that are ``separated'' from
the origin, in the sense that $D\cap(-\varepsilon,\varepsilon
)=\varnothing$ for some $\varepsilon>0$. If the whole path of the process
were available
(and, hence, the jumps of the process would be {observable}),
the problem would be identical to the estimation of the
intensity of a nonhomogeneous Poisson process
on a fixed time interval, say $[0,1]$,
based on $[T]$ independent copies of the process.
Unfortunately, under discrete-sampling, the times and magnitudes of jumps
are latent (unobservable) variables. {Nevertheless, it is expected
that the statistical property of the jumps can be inferred when the
frequency and time horizon of
observations increase to infinity, which is precisely the sampling
scheme we adopt in this paper}.

Nonparametric estimators for the L\'evy density were proposed in \cite
{FigHou2006},
under continuous sampling of the process, and in \cite{Fig2007a},
under discrete sampling,
using the \textit{method of sieves}.
The method of sieves was originally proposed by Grenander \cite{Grenander}
and has been applied more recently by Birg\'e, Massart and others
(see, e.g., \cite{Birge1999,Birge})
to several classical nonparametric problems, such as
density estimation and regression.
This approach consists of the following general steps.
First, choose a family of finite-dimensional \textit{linear models}
of functions, called \textit{sieves}, with good approximation
properties. Common sieves are splines, trigonometric polynomials
and wavelets.
Second, specify a ``distance'' metric $d$ between functions, relative
to which the best approximation of $s$
in a given linear model $\calS$
will be characterized. That is, the best approximation $s^{\bot}$ of
$s$ on $\calS$ is given by
$
d(s,s^{\bot})=\inf_{p\in\calS} d(s,p).
$
Finally, devise an estimator $\hat{s}$,
called the \textit{projection estimator},
for the best approximation $s^{\bot}$ of $s$ in
$\calS$.

The sieves considered here are of the general form
%
\begin{equation}\label{ApprxmModel}
\calS:=\{\beta_{1}\varphi_{1}+\cdots+\beta_{d}\varphi_{d}\dvtx
\beta_{1},\ldots,\beta_{d}\in\bbr\},
\end{equation}
where $\varphi_{1},\ldots, \varphi_{d}$ are
orthonormal functions with respect to the inner product
$
\langle p,q\rangle _{D}:=\int_{D} p(x) q(x)\, \mathrm{d}x.
$
In the sequel,
$\|\cdot\|:=\|\cdot\|_{D}$ stands for the associated
norm $\langle\cdot,\cdot\rangle_{D}^{1/2}$
on $\bbl^{2}(D,\mathrm{d}x)$.
We recall that, relative to the distance induced by $\|\cdot\|$,
the element of $\calS$ closest to $s$,
that is, the \textit{orthogonal projection} of $s$ on $\calS$,
is given by
%
\begin{equation}\label{OrthgProj}
s^{\bot}(x):=\sum_{j=1}^{d}\beta(\varphi_{j})\varphi_{j}(x),
\end{equation}
where
$\beta(\varphi_{j}):=\langle \varphi_{j},s\rangle _{D}=\int_{D} \varphi
_{j}(x) s(x) \,\mathrm{d}x$.
Thus, under this setting, the method of sieves reduces to the
estimation of the functional
\[
\beta(\varphi)=\int_{D} \varphi(x) s(x) \,\mathrm{d}x
\]
for certain functions $\varphi$. In Section \ref{SectCLT}, we propose
estimators for $\beta(\varphi)$ and, as a by-product, we develop
projection estimators $\hat{s}$ on $\calS$.

Following \cite{Fig2007a}, we further specialize our approach and take
\textit{regular piecewise polynomials} as sieves, although similar
results will hold true if we take other typical classes of sieves, such
as smooth splines, trigonometric polynomials or wavelets.
For future reference, let us formally define the sieves.
\begin{defn}\label{RegPolSiev}
$\calS_{k,m}$ stands for the class of functions $\varphi$ such that
for each $i=0,\ldots, m-1$,
there exists a polynomial $q_{i,k}$ of degree at most $k$ such that
$\varphi(x)=q_{i,m}(x)$ for
all $x$ in $(x_{i-1},x_{i}]$, where $x_{i}=a+i(b-a)/m $.
\end{defn}

It is easy to build an orthonormal basis for $\calS_{k,m}$ using the
orthonormal Legendre polynomials $\{Q_{j}\}_{j\geq 0}$ on $\bbl
^{2}([-1,1],\mathrm{d}x)$. Indeed, the functions
%
\begin{equation}\label{BasisSplines}
\hat{\varphi}_{i,j}(x):=
\sqrt{\frac{2j+1}{x_{i}-x_{i-1}}}
Q_{j}\biggl(\frac{2x-(x_{i}+x_{i-1})}{x_{i}-x_{i-1}}\biggr)
\mathbf{1}_{[x_{i-1},x_{i})}(x),
\end{equation}
with $i=1,\ldots,m$ and $j=0,\ldots,k$, form an orthonormal basis for
$\calS_{k,m}$. For future reference, let us recall that
%
\begin{equation}\label{EstimateLegendre}
|Q_{j}(x)|\leq 1
\quad \mbox{and}\quad
|Q_{j}'(x)|\leq Q_{j}'(1)=\frac{j(j+1)}{2}.
\end{equation}

We now review a few points of \cite{Fig2007a} in order to motivate the
results in this paper.
It is proved in \cite{Fig2007a} that by appropriately choosing the
number of classes $m$ and the sampling frequency high enough (both
choices determined as a function of the time horizon $T$), the resulting
projection estimator on $\calS_{m,k}$ attains the same rate of
convergence in $T$ as the minimax risk on a certain class $\Theta$ of
smooth functions.
Specifically, the referred minimax risk, defined by
%
\begin{equation}\label{PhM}
\inf_{\hat{s}_{T}}\sup_{{s}\in\Theta}
\bbe_{s} \biggl[\int_{a}^{b}\bigl(\hat{s}_{T}(x)-s(x)
\bigr)^{2}\,\mathrm{d}x\biggr],
\end{equation}
where the infimum is over all estimators $\hat{s}_{T}$ based on $\{
X_{t}\}_{t\leq T}$, converges to $0$ at a rate $\mathrm{O}(T^{-2\alpha/(2\alpha
+1)})$ as $T\to\infty$ (see \cite{Fig2007a}, Theorem 4.2). The
parameter $\alpha$ characterizes the smoothness of the L\'evy densities
$s\in\Theta$ on the interval $[a,b]$, in that if
$s$ is $r$-times differentiable on $(a,b)$ ($r=0,\ldots$) and
%
\begin{equation}\label{HolderCond}
\bigl|s^{(r)}(x)-s^{(r)}(y)\bigr|\leq  L |x-y|^{\kappa}
\end{equation}
for all $x,y\in(a,b)$ and some $L<\infty$ and
$\kappa\in(0,1]$, then the smoothness parameter of $s$ is
$\alpha:=r+\kappa$.
In \cite{Fig2007a}, Proposition 3.5, we show that there exists a
critical mesh $\delta_{T}>0$ such that
if the time span between consecutive sampling observations is at most
$\delta_{T}$ and $m_{T}:=[T^{1/(2\alpha+1)}]$, then
the resulting projection estimator, denoted by $\widetilde{s}_{T}$,
is such that
\begin{equation}
\label{RateProjEst}
\limsup_{T\rightarrow \infty}
T^{{2\alpha}/{(2\alpha+1)}}
\sup_{{s}\in\Theta}
\bbe\|s - \widetilde{s}_{T}\|^{2}<\infty.
\end{equation}
Of course, an ``explicit'' estimate of $\delta_{T}$ is necessary for
practical reasons. In Section \ref{CriticalMeshSect}, we show that it
is sufficient that $\delta_{T}=\mathrm{O}(T^{-1})$, improving a former result
in \cite{Fig2007a} (see Proposition 3.7 therein).

Note that the convergence in (\ref{RateProjEst}) is in the integrated
mean square sense.
A natural question, one which we consider in this paper, is whether or
not projection estimators $\hat{s}_{T}$ on
$\calS_{k,m}$ can be devised such that
%
\begin{equation}\label{ExactCLT}
T^{\alpha/(2\alpha+1)}\bigl(\hat{s}_{T}(x)-s(x)\bigr)
\,\ld\,\bar{\sigma}(x) Z
\end{equation}
holds for a standard normal random variable $Z$, for each fixed $x\in D$.
We were unable to obtain~(\ref{ExactCLT}) due to the fact that
the bias of the estimator $\hat{s}_{T}$,
namely $\bbe\hat{s}_{T}(x) -s(x)$, is
just $\mathrm{O}(T^{-\alpha/(2\alpha+1)})$. However, for any
$\beta<\frac{\alpha}{2\alpha+1}$, we can devise a projection estimator
$\hat{s}_{T}^{\beta}$ such that\vspace{-5pt}
%
\begin{equation}\label{WeakCLT}
T^{\beta}\bigl(\hat{s}_{T}^{\beta}(x)-s(x)\bigr)
\,\ld\,\bar{\sigma}(x) Z.
\end{equation}
The idea is to use ``undersmoothing'' to make the effect of bias
negligible. Our results are in keeping with those obtained in other
standard nonparametric
problems, such as density estimation and
functional regression, using local nonparametric methods such as kernel
estimation
(see, e.g., \cite{Hall92}). We were unable to find a reference
where undersmoothing is used in a global nonparametric method such as
the sieves method and, hence, this could be an additional contribution
of the results presented here.

An important extension of the pointwise central limit theorems is the
development of global measures of deviation or asymptotic confidence
bands for the L\'evy density. In this paper, we establish these methods
for piecewise constant and piecewise linear regular polynomials
(although we believe the result holds true for a general degree),
following ideas of the seminal work of Bickel and Rosenblatt \cite
{BicRos}. There are some important differences, however, starting from
the fact that Bickel and Rosenblatt
considered kernel estimators for probability densities, while, here, we
consider a global nonparametric method. In spite of these differences,
our results are consistent with the case of density estimation,
achieving a convergence rate of order arbitrarily close to $\log^{-1/2}(n)
\cdot n^{-1/3}$, where $n$ is the number of observations. Again, the
rate is valid
provided that the time horizon $T_{n}$ and the dimension of the sieves
is appropriately
chosen.

The paper is structured as follows. In Section \ref{CriticalMeshSect},
we derive a short-term ergodic property of a L\'evy process, which
plays a fundamental role in our results.
In Section \ref{SectCLT}, we introduce the projection estimators for
the L\'evy densities and show pointwise central limit theorems for them.
The uniform case and the resulting confidence bands are developed in
Section \ref{SectCBnds}.
Section \ref{VGM} illustrates the performance of the projection
estimators and confidence bands using a simulation experiment in the
case of a variance gamma L\'evy model. Finally, two appendices collect
the technical details of our results.

\section{An useful small-time asymptotic result}\label{CriticalMeshSect}
The critical time span $\delta_{T}$ required for the validity of
(\ref{RateProjEst}) was characterized
in \cite{Fig2007a} by the property that\vspace{-5pt}
%
\begin{equation}\label{CrtclMesh}
\sup_{y\in D} \biggl|
\frac{1}{\Delta}\bbp[X_{\Delta}\geq y]-\nu([y,\infty
))\biggr|
< k  \frac{1}{T}
\end{equation}
for all $0<\Delta<\delta_{T}$, where $k$ is a constant
(independent of $T$ and $\Delta$).
For practical reasons, an ``explicit'' estimate of this critical mesh
is necessary.
The following proposition shows that $\delta_{T}=T^{-1}$ suffices
and serves as the fundamental property of L\'evy processes used for the
asymptotic theory developed in this paper.
The proof of the proposition is provided in Appendix \hyperref[SectAppdxA]{A};
also, see~\cite{FigHou2008} for related higher order polynomial
expansions for $\bbp(X_{t}\geq y)$.
\begin{prop}\label{UnfrmCvrgnc}
Suppose that the L\'evy density $s$ of $X$
is
{Lipschitz in} an open set $D_{0}$ containing
$D=[a,b]\subset\bbr\backslash\{0\}$ and that $s(x)$ is uniformly bounded
on $|x|>\delta$ for any $\delta>0$.
Then, there exist a $k>0$ and a $t_{0}>0$ such that,
for all $0<t<t_{0,}$
%
\begin{equation}\label{UnifBound}
\sup_{y\in D} \biggl|
\frac{1}{t}\bbp[X_{t}\geq y]-
\nu([y,\infty))\biggr|
< k  t.
\end{equation}
\end{prop}

\section{Pointwise central limit theorem}\label{SectCLT}
Throughout this paper, we assume that the L\'evy process $\{X_{t}\}
_{t\geq 0}$
is being sampled over a time horizon $[0,T]$
at discrete times $0=t_{T}^{0}<\cdots<t^{n_{T}}_{T}=T$.
We also use the notation $\pi_{T}:=\{t^{k}_{T}\}
_{k=0}^{n_{T}}$ and
$\bar{\pi}_{T}:=\max_{k} \{t^{k}_{T}-t^{k-1}_{T}\}$,
where we will sometimes drop the subscript $T$.
The following statistics are the main building blocks for our estimation:
%
\begin{equation}\label{DscrtStat1}
\hat\beta^{\pi_{T}}(\varphi):=\frac{1}{T}
\sum_{k=1}^{n_{T}} \varphi
(X_{t^{k}_{T}}-X_{t^{k-1}_{T}}).
\end{equation}
In the case of a quadratic function $\varphi(x)=x^{2}$,
$ \sum_{k=1}^{n_{T}} \varphi
(X_{t^{k}_{T}}-X_{t^{k-1}_{T}})$
is the so-called realized quadratic variation of the process.
Thus, the statistics (\ref{DscrtStat1}) can be interpreted as the
realized $\varphi$-variation of the process per unit time based on the
observations $X_{t^{0}_{T}},\ldots,X_{t^{n_{T}}_{T}}$.
The estimators~(\ref{DscrtStat1}) were proposed independently by
Woerner \cite{Woerner2003} and Figueroa-L\'opez~\cite{Figueroa2004}.

The main virtue of the statistics (\ref{DscrtStat1}) lies in its
application to recover $\beta(\varphi):=\int\varphi(x) s(x)\, \mathrm{d}x$ as $T\to
\infty$ and $\bar\pi_{T}\to 0$ for bounded $\nu$-continuous
functions $\varphi$ such that $\varphi(x)\to 0$ fast enough as $x\to
0$. This result was obtained in \cite{Woerner2003} (Theorem 5.1
therein) for regular sampling schemes and in \cite{Fig2009b}
(Proposition 2.2 therein) for general sampling schemes and a more
general class of functions $\varphi$ (see also \cite{Fig2007a},
Theorem 2.3, for related central limit theorems).
The consistency of $\hat\beta^{\pi}(\varphi)$ for $\beta(\varphi)$
leads us to propose
%
\begin{equation}\label{ProjEstm}
\hat{s}^{\pi}(x):=\sum_{j=1}^{d}\hat\beta^{\pi}(\varphi_{j})\varphi_{j}(x)
\end{equation}
as a natural estimator for the orthogonal projection $s^{\bot}$
defined in (\ref{OrthgProj}). The nonparametric estimator (\ref
{ProjEstm}) was proposed in \cite{Figueroa2004}, where the problem of
model selection was also considered under continuous-time sampling.

As was discussed in the \hyperref[sec1]{Introduction}, one can construct a projection
estimator $\widetilde{s}_{T}$
on the regular piecewise polynomials $\calS=\calS_{k,m}$ of Definition
\ref{RegPolSiev} that converges to $s$,
under the integrated mean square distance,
at a rate at least as good as $T^{-2\alpha/(2\alpha+1)}$. Such a rate
can be
ensured by ``tuning'' the number of classes $m$ in the
sieve, as
well as the sampling frequency $\bar{\pi}$,
to both the degree of smoothness $\alpha$ of $s$ and
the time horizon $T$.
It is natural to wonder whether it is possible to construct a
projection estimator
$\hat{s}_{T}$ such that
\[
T^{\alpha/(2\alpha+1)}\bigl(\hat{s}_{T}(x)-s(x)\bigr)
\,\ld\,\bar{\sigma} Z
\]
as $T\rightarrow\infty$, for $Z\sim\mathcal{N}(0,1)$
and a constant $\bar{\sigma}$.
We are unable to obtain this result due to the fact that the bias $\bbe
\hat{s}_{T}(x) -s(x)$ of any projection estimator
$\hat{s}_{T}$ is, at best, $\mathrm{O}(T^{-\alpha/(2\alpha+1)})$. However, in
this section, we show that for any
$0<\beta<\frac{\alpha}{2\alpha+1}$, there exists a projection estimator
$\hat{s}_{T}^{\beta}$ such that
\[
c'_{T}\bigl(\hat{s}_{T}^{\beta}(x)-s(x)\bigr)
\,\ld\,\bar{\sigma} Z
\]
for a normalizing constant $c'_{T}\asymp T^{\beta}$ (i.e.,
$\underline{k}  T^{\beta}\leq  c'_{T}\leq\bar{k}  T^{\beta}$ for
some constants $\underline{k},\bar{k}\in(0,\infty)$ independent of
$T$). As it is often the case, our approach consists of first obtaining a central
limit theorem for $\hat{s}(x)$ centered at $\bbe\hat{s}(x)$ with
normalizing constants $c'_{T}\asymp T^{\beta}$ and, subsequently,
making the bias $\bbe\hat{s}(x)-s(x)$ to be $\mathrm{o}(c_{T}^{-1})$. The
central limit theorem for $\hat{s}(x)$ follows from a classical central
limit theorem for row-wise independent arrays.

Below, Legendre polynomials $\{Q_{j}\}_{j\geq 0}$ on $\bbl
^{2}([-1,1],\mathrm{d}x)$ are used to devise an orthonormal basis for the sieve
$\mathcal{S}_{k,m}$ of Definition \ref{RegPolSiev}.
Also, we consider L\'evy densities $s$ whose restrictions to $D:=[a,b]$
belong to the Besov class $\calB_{\infty}^{\alpha}(L^{\infty}([a,b]))$
(i.e., functions satisfying (\ref{HolderCond})
with $r\in\mathbb{N}$ and $\kappa\in(0,1]$ such that $\alpha=r+\kappa
$). The following is the main theorem of this section. Its proof is
deferred to Appendix \hyperref[SectAppdxB]{B}.

\begin{thrm}\label{TPWC}
Suppose that the L\'evy density $s$ of $X$ satisfies the conditions of
Proposition \ref{UnfrmCvrgnc} and belongs to $\mathcal{B}^{\alpha
}_{\infty}(L^{\infty}([a,b]))$ for some $\alpha\geq 1$. Let
$c_{T}$ be a normalizing constant and let $\hat{s}_{T}$ be the
projection estimator on $\mathcal{S}_{k,m_{T}}$ based on sampling
times $\pi_{T}$ such that the following conditions are satisfied:
\begin{eqnarray*}
\mathrm{(i)}\quad
c_{T}&\stackrel{T\rightarrow\infty}{\longrightarrow}&
 \infty;\qquad
\mathrm{(ii)}\quad
\frac{c^{2}_{T}m_{T}}{T}\stackrel{T\rightarrow\infty
}{\longrightarrow}
 1;\qquad
\mathrm{(iii)} \quad  c_{T}m_{T}\bar\pi_{T}\stackrel{T\rightarrow\infty
}{\longrightarrow}0;\\
\mathrm{(iv)}\quad
c_{T} m_{T}^{-\alpha}
&\stackrel{T\rightarrow\infty}{\longrightarrow}&
 0;\qquad
\mathrm{(v)}\quad  k\geq\alpha-1.
\end{eqnarray*}
Then, for any fixed $x\in(a,b)$ for which $s(x)>0$,
%
\begin{equation}\label{MPWCLT}
\frac{c_{T}}{b_{k,m_{T}}(x)}
\bigl(\hat{s}_{T}(x)-
s(x)\bigr)
\,\ld\,\bar{\sigma}(x) Z,
\end{equation}
where
\begin{eqnarray*}
Z&\sim&\mathcal{N}(0,1),\qquad \bar\sigma^{2}(x):=(b-a)^{-1}s(x), \\
b_{k,m}^{2}(x)&:=&\sum_{j=0}^{k}{(2j+1)}
\sum_{i=1}^{m}Q_{j}^{2}\biggl(
\frac{2x-(x_{i}+x_{i-1})}{x_{i}-x_{i-1}}\biggr){\mathbf{1}_{[x_{i-1},x_{i})}(x)}.
\end{eqnarray*}
Also, for any
fixed $0<\beta<\frac{\alpha}{2\alpha+1}$, the resulting projection
estimator $\hat{s}_{T}$ with
$m_{T}=[T^{1-2\beta}]$ is such that
\[
\frac{T^{\beta}}{b_{k,m_{T}}(x)}\bigl(\hat{s}_{T}(x)-s(x)\bigr)
\,\ld\,\bar\sigma(x) Z,
\]
provided that
$\bar\pi_{T}=T^{-\gamma}$ with $\gamma>1-\beta$.
\end{thrm}

\begin{remk}
\begin{enumerate}[(1)]
\item[(1)]
In view of (\ref{EstimateLegendre}), $1\leq b_{k,m} \leq\sum
_{j=0}^{k} (2j+1)$ and, hence, the normalizing constant
$c'_{T}:=c_{T}/b_{k,m_{T}}\asymp c_{T}$. Also, note that
$b_{k,m}\equiv1$ in the piecewise constant case ($k=0$).
\item[(2)] Theorem \ref{TPWC} will allow us to construct approximate
confidence intervals for $s(x)$. Concretely, the $100(1-\alpha)\%$
interval for $s(x)$ is approximately given by
\[
\hat{s}_{T}(x) \pm\frac{b_{k,m_{T}}(x)}{c_{T}(b-a)^{1/2}}\hat
{s}_{T}^{1/2}(x) z_{\alpha/2},
\]
where $z_{\alpha/2}$ is the $\alpha/2$ normal quantile.
\end{enumerate}
\end{remk}

\section{Confidence bands for L\'evy densities}\label{SectCBnds}
In this section, we address the problem of constructing confidence
bands for the L\'evy density $s$ of a L\'evy process using projection
estimators $\hat{s}_{T}^{n}$ on $\calS_{k,m}$ based on $n$
evenly-spaced observations~of the process at $t_{0}=0<\cdots<t_{n}=T$ on
$[0,T]$. Confidence bands entail the limit in distribution of the
uniform norm
\[
\|\hat{s}^{n}_{T} -s\|_{[a,b]}:=\sup_{x\in[a,b]} |\hat
{s}^{n}_{T}(x) -s(x)|,
\]
but, as before, we will first work with the uniform norm of
%
\begin{equation}\label{DP}
Y^{n}_{T}(x)
:=\hat{s}^{n}_{T}(x)-\bbe\hat{s}^{n}_{T}(x),\qquad  x\in[a,b],
\end{equation}
and then estimate the uniform norm of the bias $\bbe\hat
{s}^{n}_{T}(x)-s(x)$.
We follow ideas from the seminal paper of Bickel and Rosenblatt \cite
{BicRos}, wherein confidence bands for probability densities are
constructed based on kernel estimators. There are two fundamental
general directions in Bickel and Rosenblatt's approach:
\begin{enumerate}[(1)]
\item[(1)] the statistics of interest are expressed in terms of the
so-called uniform standardized empirical process
%
\begin{equation}\label{USEP}
Z_{n}^{0}(x):=n^{1/2}\{F_{n}^{*}(x)-x\}, \qquad x\in[0,1],
\end{equation}
where, denoting by $F_{t}$ the distribution of $X_{t}$ and by $\delta
^{n}:=t_{i}-t_{i-1}$ the time span between observations,
$F_{n}^{*}(\cdot)$ is the empirical distribution of $\{F_{\delta
^{n}}(X_{t_{i}}-X_{t_{i-1}})\}_{i\leq n}$;
\item[(2)] the empirical process $Z_{n}^{0}$ is approximated by a
Brownian bridge $Z^{0}$ and the error is estimated using Brillinger's
result \cite{Bril} or the Koml\'os, Major and Tusn\'ady construction
\cite{KMT}.
\end{enumerate}

Once the statistic of interest is related to the Brownian bridge
$Z^{0}$, we will carry over several successive approximations (see
Appendix  \hyperref[SectAppdxC]{C} for the details), which will allow the
distribution of $\|Y^{n}_{T}\|_{[a,b]}$ to be connected with the
limiting distribution of the extreme value
\[
\bar{M}_{m}:=\max_{1\leq j\leq m}\bigl\{\zeta^{(k)}_{j}\bigr\}
\]
of independent copies $\{\zeta_{j}^{(k)}\}_{j}$ of the random variable
%
\begin{equation}\label{EqKeySprm}
\zeta^{(k)}:=\sup_{x\in[-1,1]}\Biggl|\sum_{j=0}^{k}{\sqrt{2j+1}}
Q_{j}(x) Z_{j}\Biggr|,
\end{equation}
where $Z_{j}$ are i.i.d. standard normal random variables. The problem
is then reduced to finding the extreme value distribution of a random
sample from (\ref{EqKeySprm}). For instance, in the case $k=0$, $\zeta
^{(0)}_{j}\iid|Z_{0}|$, which is known to satisfy
%
\begin{equation}\label{BscLmt}
\lim_{n\rightarrow\infty}\bbp\biggl(\max_{1\leq j\leq m} \bigl|\zeta
_{j}^{(0)}\bigr| \leq\frac{{y}}{a_{m}}+b_{m}\biggr)=\mathrm{e}^{-2\mathrm{e}^{-{y}}}
\end{equation}
for any $y>0$, where
\begin{eqnarray}\label{NrmlzCnst1}
a_{m}&=&(2\log m)^{1/2},\\
\label{NrmlzCnst2}
b_{m}&=&(2\log m)^{1/2}-\tfrac{1}{2}(2\log m)^{-1/2}
(\log\log m +\log4 \uppi).
\end{eqnarray}
We are also able to tackle the case $k=1$, where $\zeta
^{(1)}=|Z_{0}|+\sqrt{3}|Z_{1}|$, but the general case is still under
investigation. Our assumptions are as follows.
\begin{cond}\label{CndLevyDnsty}
\begin{enumerate}[(1)]
\item[(1)]$s$ is positive and continuous on $[a,b]$.
\item[(2)]$s$ is differentiable in $(a,b)$ and, moreover, the derivative of
$s^{1/2}$ is bounded in absolute value on $(a,b)$.
\end{enumerate}
\end{cond}

We are ready to present the main result of this section. We defer its
proof to Appendix~\hyperref[SectAppdxC]{C}.
\begin{thrm}\label{TCB}
Suppose that $\nu(\bbr)=\infty$ or $\sigma\neq0$. Also, suppose that
the L\'evy density $s$ satisfies the conditions of Proposition \ref
{UnfrmCvrgnc} and the Assumption~\textup{\ref{CndLevyDnsty}}. Let
$T_{n}\to \infty$ and $m_{n}\to\infty$ be such that
\[
\textup{(i)}\quad
\delta^{n}\log\delta^{n} \cdot m_{n}\log m_{n}\stackrel{n\rightarrow
\infty}{\longrightarrow}
 0,\qquad
\textup{(ii)} \quad \frac{\log^{2}n}{T_{n}} \cdot m_{n}\log m_{n}\stackrel
{n\rightarrow\infty}{\longrightarrow}
 0,
\]
where $\delta_{n}:=T_{n}/n$.
Then, for $k\in\{0,1\}$, the deviation process $Y_{T_{n}}^{n}$ of
(\ref{DP}) satisfies
%
\begin{equation}
\lim_{n\rightarrow\infty}\bbp\Bigl(
a_{m_{n}}\Bigl\{
\kappa  \bar{T}_{n}^{1/2}\sup_{x\in[a,b]} |
s^{-1/2}(x)Y_{T_{n}}^{n}(x)|-b_{m_{n}}\Bigr\}\leq{y}
\Bigr)=\mathrm{e}^{-\kappa' \mathrm{e}^{-{y}}},
\end{equation}
where $\bar{T}_{n}:=T_{n}/m_{n}$, $a_{m}$ and $b_{m}$ are defined as in
\textup{(\ref{NrmlzCnst1})--(\ref{NrmlzCnst2})} and $(\kappa,\kappa
')=((b-a)^{1/2},2)$ if $k=0$ or $(\kappa,\kappa
')=((b-a)^{1/2}2^{-1},4)$ if $k=1$.
\end{thrm}

The previous result shows that
\[
a_{m_{n}}\Bigl\{\kappa  \bar{T}_{n}^{1/2}\sup_{x\in[a,b]}
s^{-1/2}(x) | \hat{s}^{n}_{T_{n}}(x)-\bbe\hat
{s}^{n}_{T}(x)|-b_{m_{n}}\Bigr\}
\]
converges to a Gumbel distribution. The final step in constructing our
confidence bands consists of finding conditions for replacing $\bbe\hat
{s}^{n}_{T}$ with $s$. The following result shows this step. Its
proof is presented in Appendix  \hyperref[SectAppdxC]{C}.

\begin{cllry}\label{KRCB}
Suppose that the conditions of Theorem \ref{TCB} hold true,
that the restriction of $s$ to $[a,b]$ is a member of $\mathcal
{B}^{\alpha}_{\infty}(L^{\infty}([a,b]))$ and also that
%
\begin{equation}\label{RAlph0}
\mathrm{(iii)}\quad   T_{n} m_{n}^{1-2\alpha} \log^{2} m_{n}\stackrel{n\rightarrow
\infty}{\longrightarrow}
 0.
\end{equation}
Then,
%
\begin{equation}\label{EqLimitMaxima3}
{\lim_{n\rightarrow\infty}\bbp\biggl(
a_{m_{n}}\biggl\{\kappa  \bar{T}_{n}^{1/2}\sup_{x\in[a,b]} \frac
{1}{s^{1/2}(x)}| \hat{s}^{n}_{T_{n}}(x)-s(x)
|-b_{m_{n}}\biggr\}
\leq{y}\biggr)=\mathrm{e}^{-\kappa' \mathrm{e}^{-{y}}}},
\end{equation}
where we have used the same notation {for $\kappa$ and $\kappa'$}
as in Theorem \ref{TCB}.
\end{cllry}

The previous corollary allows us to construct confidence bands for $s$
on $[a,b]$ based on the projection estimators $\hat{s}$ on regular
piecewise linear (or constant) polynomials.
Indeed, suppose that $y_{\alpha}^{*}$ is such that
$
\exp\{-k'\mathrm{e}^{-y_{\alpha}^{*}}\}=1-\alpha
$
and let
\[
{d_{n}:=\frac{1}{\sqrt{2}\kappa}\biggl(\frac{y^{*}_{\alpha
}}{a_{m_{n}}}+b_{m_{n}}\biggr)
\bar{T}_{n}^{-1/2}}.
\]
Then, as $n\rightarrow\infty$,
%
\begin{equation}\label{EqCI}
s(x)\in
\bigl(\hat{s}^{n}_{T_{n}}(x)+\bigl\{d_{n}^{2}\pm
\sqrt{\bigl(\hat{s}^{n}_{T_{n}}(x)+d_{n}^{2}\bigr)^{2}-(\hat
{s}^{n}_{T_{n}}(x))^{2}}\bigr\}\bigr),
\end{equation}
with $100(1-\alpha)\%$ confidence. The above interval is asymptotically
equivalent to the following, simpler, interval:
%
\begin{equation}\label{EqCIb}
{s(x)\in
\biggl(\hat{s}^{n}_{T_{n}}(x)\pm
\frac{1}{\kappa}\biggl(\frac{y^{*}_{\alpha}}{a_{m_{n}}}+b_{m_{n}}
\biggr)\bar{T}_{n}^{-1/2}
(\hat{s}^{n}_{T_{n}}(x))^{1/2}\biggr)}.
\end{equation}

We conclude this section with some final remarks.
\begin{remk}
In the case where $T_{n}:= c_{n} \cdot n^{\alpha_{1}}$ and $m_{n}=
[ d_{n} \cdot n^{\alpha_{2}}]$,
for some $\alpha_{1},\alpha_{2}>0$, $c_{n}\asymp1$ and $d_{n} \asymp
1$, the conditions \textup{(i)--(ii)} of Theorem \ref{TCB} are satisfied if
$0<\alpha_{1}<1$ and $0<\alpha_{2}<(1-\alpha_{1})\wedge\alpha_{1}$.
Also, it can be checked that condition \textup{(iii)} of Corollary \ref{KRCB} is
met if
%
\begin{equation}\label{RAlph}
{0<\alpha_{1}<
\frac{2\alpha+1}{3\alpha+2}\quad   \mbox{and}\quad
\frac{\alpha_{1}}{1+2\alpha}<\alpha_{2}<(2-3\alpha_{1})\wedge\alpha_{1}}.
\end{equation}
Note that $(\alpha_{2}-\alpha_{1})/2$ can be made arbitrarily close to
$-\alpha/(3\alpha+1)$ on the range of values \textup{(\ref{RAlph})} and, thus,
$a_{m_{n}} \bar{T}_{n}^{-1/2}$ can be made to vanish at a rate
arbitrarily close to
$
( \log n)^{-1/2} n^{-{\alpha}/{(3\alpha+1)}},
$
provided that $\alpha$ is large enough.
In particular, if $0<\varepsilon\ll1$ and $s$ is smooth enough, then
$m_{n}$ and $T_{n}$ can be chosen such that
\[
\|\hat{s}_{T_{n}}^{n}-s\|_{[a,b]}=\mathrm{O}(\log^{-1/2}(n)  n
^{-1/3+\varepsilon}).
\]
\end{remk}

\section{A numerical example}\label{VGM}

Variance gamma processes (VG) were proposed in {\cite
{Madan3} and \cite{Madan1} as substitutes for
Brownian motion in the Black--Scholes model. Since their
introduction, VG processes have received a great
dealt of attention, even in the financial industry.
A variance gamma
process $X=\{ X(t)\}_{t \geq0}$ is a time-changed Brownian motion with
drift of the form
%
\begin{equation}
\label{TimeChangeDef}
X(t)= \theta U(t) + \sigma W(U(t)),
\end{equation}
where $\{ W(t)\}_{t\geq0}$ is a standard Brownian motion, $\theta\in
\mathbb{R}$, $\sigma>0$ and $U=\{U(t)\}_{t \geq0}$ is an independent
gamma L\'evy process such that
$\mathrm{E}[U(t)]=t$ and $\operatorname{Var}[U(t)]=\nu t$.
Since gamma processes are \textit{subordinators}, the process $X$ is
itself a L\'evy process (see \cite{Sato}, Theorem 30.1) and its
L\'evy density takes the form
%
\begin{equation}\label{LevyDnstyVG}
s(x)=\cases{
\displaystyle\frac{\alpha}{|x|}\exp\biggl(-\frac{|x|}{\beta^{-}}\biggr), & \quad $\mbox{if }
x<0,$\vspace*{2pt}
\cr
\displaystyle\frac{\alpha}{x} \exp\biggl(-\frac{x}{\beta^{+}}\biggr), & \quad $\mbox{if }
x>0,$}
\end{equation}
where $\alpha>0$, $\beta^{-} \geq0$ and $\beta^{+} \geq0$ with
$|\beta^{-}|+|\beta^{+}|>0$ (see, e.g., \cite{Cont2003} for
expressions for $\beta_{\pm},\alpha$ in terms of $\theta$, $\sigma$ and
$\nu$).
In that case,
$\alpha$ controls the overall jump activity, while
$\beta^{+}$ and $\beta^{-}$ take charge of
the intensity of large positive and negative jumps, respectively. In particular,
the difference between $1/\beta^{+}$ and $1/\beta^{-}$ determines the
frequency of drops relative to rises, while their sum measures the
frequency of
large moves relative to small ones.

The performance of projection estimation for the variance gamma L\'evy
process was illustrated in \cite{Fig2007a} via simulation experiments.
In this section, we want to further extend this analysis to show
the performance of confidence bands. As in \cite{Fig2007a},
we take as sieve the class $\calS_{0,m}$, namely, the span of the
indicator functions $\chi_{[x_{0},x_{1}]} , \ldots,
\chi_{(x_{m-1},x_{m}]}$, where $x_{0}<\cdots<x_{m}$ is a regular
partition of an interval $D \equiv[a,b]$, with $0<a$ or $b<0$.
We take parameter values which are partially motivated by the empirical
findings of \cite{Madan1} based on daily returns of the S\&P500
index from January 1992 to September 1994
(see their Table~I).
Using maximum likelihood methods,
the annualized estimates of the parameters for
the variance gamma model were reported to be
$\hat{\theta}_{\mathit{ML}}=-0.00056256$, $\hat{\sigma}^{2}_{\mathit{ML}}=0.01373584$
and $\hat{\nu}_{\mathit{ML}}=0.002$, from which it can easily be found that
%
\begin{equation}\label{VGP}
\hat\alpha= 500, \qquad \hat{\beta}^{+}=0.0037056 \quad \mbox{and}\quad
\hat{\beta}^{-}=0.0037067.
\end{equation}
These parameter values seem to be consistent with other empirical
studies (see, e.g., \cite{Seneta}), although we admit that
parameter values fitted to intraday high-frequency data would have been
preferable.

We simulate 100 samples of the VG process with a maximal time horizon
of $T=10$ years and a sampling span between observations of $\delta
=1/(252\times 6.5\times 60\times 12)$. Assuming a business calendar year of 252 days and
a trading day of 6.5 hours, the time span between observations
corresponds to $5$ seconds. Intraday data of such characteristics is
available via financial databases such as NASDAQ TAQ.

We estimate the sample coverage probabilities
\[
c_{\alpha}:=\bbp\bigl(s(\cdot)\in\mbox{the }100(1-\alpha)\%\mbox{
confidence band on }[a,b]\bigr),
\]
based on the 100 simulations for two sampling frequencies $\delta
=1/(252\times 6.5\times 60\times 12)$ (5~seconds) and
$\delta=1/(252\times 6.5\times 60)$ (1 minute), and maturities of $T=1,3,5$ and
$10$ years.
We use two possible numbers of classes: $m=40$ and the data-driven
selected $m$ proposed in \cite{Fig2007a}. Concretely, the selection
criterion is given by
\begin{equation}\label{DscrtCrtrnMdlSlctn}
\hat{m}:=\mathop{\operatorname{argmin}}_{m} \{
-\|\hat s^{\pi}_{m} \|^{2} + \operatorname{pen}^{\pi}(\calS_{k,m}) \},
\end{equation}
where $\hat{s}^{\pi}_{m}$ is given according to (\ref{ProjEstm}) and
$\operatorname{pen}^{\pi}$ is given
by
%
\begin{equation}
\label{PenaltyDiscrete1}
\operatorname{pen}^{\pi}(\calS_{k,m})=\frac{2}{T^{2}}
\sum_{i=1}^{n}
 \sum_{i,j} \hat{\varphi}_{i,j}^{2}(X_{t_{i}}-X_{t_{i-1}}).
\end{equation}
The quantity to be minimized in (\ref{DscrtCrtrnMdlSlctn}) is a
discrete-time version of an unbiased estimator of the shifted risk
$\bbe\|s-\hat{s}^{\pi}_{m}\|^{2}-\|s\|^{2}$ (see \cite{Fig2007a},
Section 5, for more details).


\begin{figure}[b]

\includegraphics{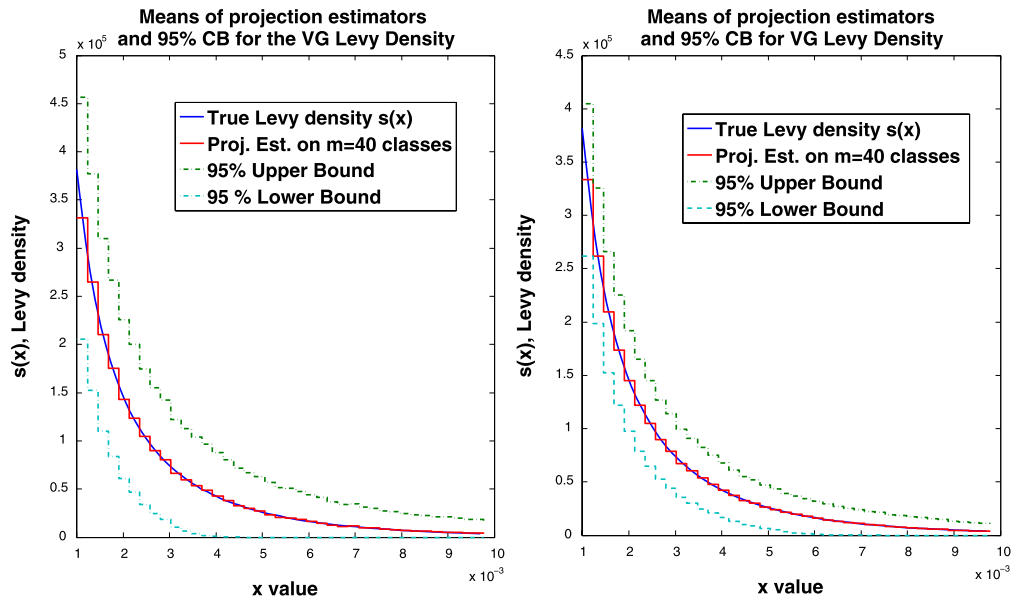}

\caption{Means of projection estimators and corresponding confidence
bands for the VG model based on $100$ simulations with a sampling time
span of $1/(252\times 6.5\times 60 \times 12)$ (about $5$ seconds) during $3$ years (left
panel) and $10$ years (right panel).}\label{CB5y5sa}
\end{figure}


The Table \ref{tab1} shows the coverage probabilities for the interval
$[a,b]=[0.001,0.1]$ (based on $100$ simulations).
Overall, the coverage probabilities of the confidence bands for $m=40$
are good. In the case of the data-driven selected $m$, there are some
values of $m$ for which probabilities are quite low. Such cases occur
(only) when the band does not contain the density very near $a=0 .001$.
It seems more reasonable to take an average between different classes
with values of~$m$ which are reasonably close in terms of the quantity
in (\ref{PenaltyDiscrete1}).
\begin{table}
  \caption{Empirical coverage probabilities of 95\% confidence bands on the interval
[$0.001,0.1$] based on a~piece-wise projection estimator with~$m$ classes}\label{tab1}
\begin{tabular*}{\textwidth}{@{\extracolsep{\fill}}lllll@{}}
\hline
\multicolumn{1}{@{}l}{$\delta\backslash T$} & 1 year & 3 years & 5
years & 10 years \\
\hline
5 s&0.97 ($m=40$)& 0.99 ($m=40$)& 0.97 ($m=40$)&   0.97 ($m=40$)\\
&0.98 ($m=35$) & 0.95 ($m=25$)&0.80 ($m=25$)&\\[3pt]
1 min& 0.93 ($m=40$)&0.94 ($m=40$)&0.98 ($m=40$)&0.87 ($m=40$)\\
&0.97 ($m=35$)& 0.75 ($m=25$)&0.60 ($m=25$)&0.94 ($m=50$) \\
\hline
\end{tabular*}
\end{table}

To illustrate how close the estimated L\'evy density is to the true L\'
evy density and the overall width of the confidence bands, Figure~\ref
{CB5y5sa} shows the actual L\'evy density (solid blue line), the mean
of the penalized projection estimator (solid red line) and the means of
the lower and upper $95 \%$-confidence bands (dashed lines). All the
means are computed using 100 confidence bands based on $\delta=5$
seconds and time horizons of $T=3$ and $T=10$ years.
The analogous figures with a sampling time span of $\delta=1$ minute
are shown in Figure~\ref{CB5y1ma}. In our empirical results (not shown
here for the sake of space), we found that high-frequency data is
crucial to estimate the L\'evy density near the origin. For instance,
the confidence bands near the origin do not perform well when taking
$30$-minute observations in a time period of $10$ years.
The Table \ref{tab2} gives the estimated coverage probabilities on the
interval $[0.005,0.2]$ based on $30$-minute returns.

Let us finish with two remarks. First, from an algorithmic point of view,
the estimation for the variance gamma model
using penalized projection is not different from
the estimation of the gamma L\'evy process.\vadjust{\goodbreak}
We can simply estimate both tails of the
variance gamma process separately.
However, from the point of view of maximum likelihood
estimation (MLE), the problem is numerically challenging.
Even though the marginal density functions have ``closed''
form expressions (see \cite{Madan1}),
there are well-documented issues with MLE
(see, e.g., \cite{Prause}).
Finally, it worth pointing out that applying an efficient estimation
method to a misspecified model could lead to quite undesirable results,
as was
illustrated in \cite{Fig2007a}, where MLE was applied to a
CGMY model (see \cite{Madan}) with parameter values quite close to
those of a gamma process.
The numerical experiments in \cite{Fig2007a} show that a modestly
efficient robust nonparametric
method is sometimes preferable to a very efficient estimation method.

\renewcommand{\thetable}{\arabic{table}}
\setcounter{table}{1}
\begin{table}
\caption{Empirical coverage probabilities of 95\% confidence bands on the interval
$[0.005,0.2]$ based on a~piece-wise projection estimator with~$m$ classes}\label{tab2}
\begin{tabular*}{\textwidth}{@{\extracolsep{\fill}}lllll@{}}
\hline
\multicolumn{1}{@{}l}{$\delta\backslash T$} & 1 year & 3 years & 5
years & 10 years \\
\hline
30 min &0.34 ($m=40$)&0.73 ($m=40$)&0.87 ($m=40$)&0.97 ($m=40$)\\
&0.43 ($m=10$)& 0.71 ($m=35$)&0.85 ($m=35$)&0.97 ($m=25$)\\
\hline
\end{tabular*}
\end{table}
%

\begin{figure}[b]

\includegraphics{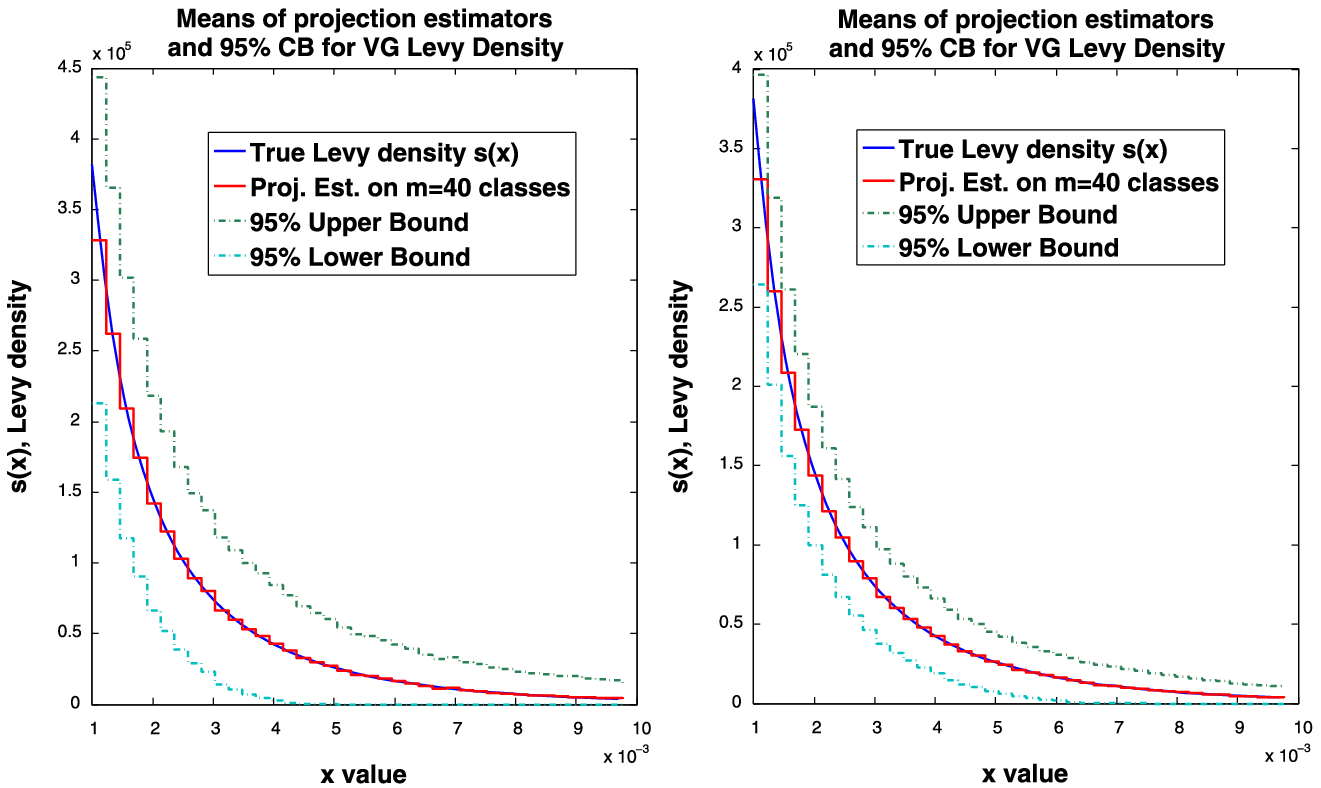}

\caption{Means of projection estimators and corresponding confidence
bands for the VG model based on $100$ simulations with a sampling time
span of $1/(252\times 6.5\times 60)$ (about $1$ minute) during $3$ years (left
panel) and $10$ years (right panel).}\label{CB5y1ma}
\end{figure}

\begin{appendix}

\section{\texorpdfstring{Proof of Proposition \protect\ref{UnfrmCvrgnc}}{Proof of Proposition 2.1}}\label{SectAppdxA}

Without loss of generality, we assume that $a>0$.
Consider the process
%
\setcounter{equation}{0}
\begin{equation}\label{TrctedLevy}
\widetilde{X}^{\varepsilon}_{t}:=
\int_{0}^{t}\!\!\int_{\bbr} x\mathbf{1}_{\{|x|\geq\varepsilon\}}
\mu(\mathrm{d}x,\mathrm{d}s)
\end{equation}
for $0<\varepsilon<1$, which is well known to be a
compound Poisson process with intensity of jumps
$\lambda_{\varepsilon}:=\nu(\{|x|\geq\varepsilon\})$
and jump distribution $\frac{1}{\lambda_{\varepsilon}}\mathbf{1}_{\{
|x|\geq\varepsilon\}}\nu(\mathrm{d}x)$.
The remainder process, $X^{\varepsilon}:=X-\widetilde{X}^{\varepsilon}$,
is then a L\'evy process with jumps bounded by
$\varepsilon$. Concretely, $X^{\varepsilon}$
has L\'evy triplet {$(\sigma^{2},b_{\varepsilon},
\mathbf{1}_{\{|x|\leq\varepsilon\}} \nu(\mathrm{d}x))$,
where $b_{\varepsilon}=b-\int_{\varepsilon<|x|\leq 1}x\nu(\mathrm{d}x)$}.
{The following tail estimate will play an important role in the sequel:
%
\begin{equation}\label{TailEstm}
\bbp(|X^{\varepsilon}_{t}|\geq z)
\leq  \exp\{\alpha z_{0}\log z_{0}\} \exp\{\alpha z-\alpha z\log
z\}
t^{z\alpha},
\end{equation}
valid for an arbitrary, but fixed, positive real $\alpha\in
(0,\varepsilon^{-1})$
and for any $t,z>0$ such that $t< z_{0}^{-1}z$, where $z_{0}$ depends
only on $\alpha$
(see \cite{Ruschendorf}, Lemma 3.2, or \cite{Sato}, Section 26,
for a proof).}

Define
\[
A_{y}(t):=
\frac{1}{t}\biggl\{
\frac{1}{t} \bbp[X_{t}\geq y]-
\nu([y,\infty))\biggr\},
\]
which, for $\varepsilon<\frac{y}{2}\wedge 1$
and after conditioning on the number of jumps,
can be written as
\begin{eqnarray*}
A_{y}(t)&=&\frac{1}{t^{2}}\bbe f_{y}(X^{\varepsilon}_{t})
\mathrm{e}^{-\lambda_{\varepsilon}t}+
\mathrm{e}^{-\lambda_{\varepsilon}t}
\int_{|x|\geq \varepsilon}\frac{1}{t}\{
\bbe f_{y}(X_{t}^{\varepsilon}+x)-f_{y}(x)\}\nu(\mathrm{d}x)\\
&&{}-\frac{1-\mathrm{e}^{-\lambda_{\varepsilon} t}}{t}
\int_{x>y}f_{y}(x)\nu(\mathrm{d}x)+
\mathrm{e}^{-\lambda_{\varepsilon} t}\sum_{n=2}^{\infty}\frac{(\lambda
_{\varepsilon})^{n}
t^{n-2}}{n!}
\bbe f_{y}\Biggl(X^{\varepsilon}_{t}+\sum_{i=1}^{n}\xi_{i}\Biggr),
\end{eqnarray*}
where $f_{y}(x)=\mathbf{1}_{x\geq y}$.
{The first term on the right-hand side of the above expression is
bounded uniformly for $y\in[a,b]$ and $t<t_{0}$, for certain
$t_{0}(\alpha)>0$, because of (\ref{TailEstm}) taking $z=a$ and $\alpha
\in(2a^{-1},\varepsilon^{-1})$.}
The last two terms in the same expression
are uniformly bounded in absolute value
by $\nu(x\geq a)$ and $\nu(|x|\geq \varepsilon)^{2}$,
respectively.
We need to show that the second term
is uniformly bounded.
Define
$
B_{y}(t):=\int_{|x|\geq \varepsilon}\{
\bbe f_{y}(X_{t}^{\varepsilon}+x)-f_{y}(x)\}\nu(\mathrm{d}x).
$
Clearly,
\begin{eqnarray*}
B_{y}(t)&:= &\int_{y-\varepsilon}^{y}
\bbp\{X_{t}^{\varepsilon}\geq y-x\} s(x)\,\mathrm{d}x
-\int_{y}^{y+\varepsilon}
\bbp\{X_{t}^{\varepsilon}< y-x\} s(x)\,\mathrm{d}x\\
&&{} +\int_{\{x<y-\varepsilon,|x|\geq\varepsilon\}}
\bbp\{X_{t}^{\varepsilon}\geq y-x\}s(x) \,\mathrm{d}x
-\int_{y+\varepsilon}^{\infty}
\bbp\{X_{t}^{\varepsilon}< y-x\} s(x)\,\mathrm{d}x.
\end{eqnarray*}
Since $s$ is bounded and integrable away from the origin,
the last two terms in the expression for $B_{y}(t)$
can be bounded in absolute value by
$
\nu\{|x|\geq\varepsilon\}
\bbp\{|X_{t}^{\varepsilon}|\geq\varepsilon\}.
$
Dividing by~$t$, this converges to $0$ in light of
{the well-known limit
%
\begin{equation}\label{FrstLimit}
\lim_{t\rightarrow 0}\frac{1}{t}
\bbp(Z_{t}\geq z)
=\nu([z,\infty)),
\end{equation}
valid for any L\'evy process $Z$ with L\'evy measure $\nu$ and any
point $z$ of continuity of $\nu$
(see, e.g., Bertoin \cite{Bertoin}, Chapter 1).}
The other two terms can be bounded
as follows:
\begin{eqnarray*}
&&\biggl|\int_{y-\varepsilon}^{y}
\bbp\{X_{t}^{\varepsilon}\geq y-x\} s(x)\,\mathrm{d}x
-\int_{y}^{y+\varepsilon}
\bbp\{X_{t}^{\varepsilon}< y-x\} s(x)\,\mathrm{d}x\biggr|
\\
&&\quad\leq K_{1}
\int_{0}^{\varepsilon}
\bbp\{|X_{t}^{\varepsilon}|\geq u\}
u \,\mathrm{d}u + K_{0} \biggl|
\int_{0}^{\varepsilon}
\bbp\{X_{t}^{\varepsilon}\geq u\} \,\mathrm{d}u
-\int_{0}^{\varepsilon}
\bbp\{X_{t}^{\varepsilon}< -u\} \,\mathrm{d}u
\biggr|,
\end{eqnarray*}
where {$K_{1}$ is the Lipschitz constant of $s$ in $D_{0}$} and $K_{0}:=
\sup_{x\in D_{0}} |s(x)|$.
Next, applying Fubini's theorem, we can write the expression
in the last line above as follows:
\[
K_{1} \tfrac{1}{2}\bbe\{(|X_{t}^{\varepsilon}|\wedge
\varepsilon)^{2}
\} +K_{0}
|\bbe h(X_{t}^{\varepsilon})|,
\]
where $h(x)=x\mathbf{1}_{|x|\leq \varepsilon} -\varepsilon\mathbf{1}_{x<-\varepsilon}+
\varepsilon\mathbf{1}_{x>\varepsilon}$.
{Using the formulas for the variance and mean of a L\'evy process, we
obtain that
\[
\sup_{0<t\leq 1}\frac{1}{t}\bbe\{(|X_{t}^{\varepsilon}|\wedge
\varepsilon)^{2}\}
\leq
\sigma^{2} +\int_{|x|\leq \varepsilon} x^{2}\nu(\mathrm{d}x)+ b_{\varepsilon
}<\infty.
\]
}
Also,
\[
\biggl|\frac{1}{t}\bbe
h(X_{t}^{\varepsilon})\biggr|\leq
\biggl|\frac{1}{t}\bbe
X_{t}^{\varepsilon}\biggr|+
\biggl|\frac{1}{t}\bbe
X_{t}^{\varepsilon} \mathbf{1}_{\{|X_{t}^{\varepsilon}|>\varepsilon\}}
\biggr|+
\varepsilon\frac{1}{t}
\bbp\{|X_{t}^{\varepsilon}|>\varepsilon\}.
\]
The last term above converges to $0$ by (\ref{TailEstm}).
The second term also vanishes since
\[
\frac{1}{t}\bigl|\bbe
X_{t}^{\varepsilon} \mathbf{1}_{\{|X_{t}^{\varepsilon}\bigr|>\varepsilon\}
}\bigr|
\leq
\biggl\{\frac{1}{t}
\bbp\{|X_{t}^{\varepsilon}|>\varepsilon\}\biggr\}^{1/2}
\biggl\{\frac{1}{t}\bbe
(X_{t}^{\varepsilon})^{2}\biggr\}^{1/2}
\rightarrow 0
\]
as $t\rightarrow 0$. Finally, using the formula for the mean of
$X_{t}^{\varepsilon}$, we have
\[
\lim_{t\rightarrow 0}\frac{1}{t}\bbe
h(X_{t}^{\varepsilon})\leq
\lim_{t\rightarrow 0}\frac{1}{t}|\bbe
X_{t}^{\varepsilon}|=|b_{\varepsilon}|.
\]
We conclude that there exists a $t_{0}$ and $K>0$
such that for $t\leq t_{0}$,
$
\sup_{y\in D}
|B_{y}(t)|/t
\leq K .
$
This completes the proof since all other terms in $A_{y}(t)$
can be easily bounded uniformly in $D$.

\section{Proofs of the pointwise central limit theorem}\label{SectAppdxB}
Throughout this section, we shall use the orthonormal basis $\{\hat
{\varphi}_{i,j}\}_{1\leq i\leq m,0\leq j\leq k}$ of (\ref
{BasisSplines}). We start our proof with following easy lemma.
\begin{lmma}\label{EstDstnc}
Suppose that $\varphi$
has support $[c,d]\subset\bbr_{+}\backslash \{0\}$,
where $\varphi$ is continuous with continuous derivative.
Then,
\[
\biggl|\frac{\bbe\varphi(X_{\Delta})}{\Delta}
-\beta(\varphi)\biggr|\leq
\biggl(|\varphi(c)|+\int_{c}^{d}|\varphi'(u)|\,\mathrm{d}u\biggr)
M_{\Delta}([c,d]),
\]
where $\beta(\varphi):=\int\varphi(x) s(x)\, \mathrm{d}x$ and $M_{\Delta
}([c,d]):=\sup_{y\in[c,d]}
|\frac{1}{\Delta}\bbp[X_{\Delta}\geq y]-\nu([y,\infty
))|$.
\end{lmma}
\begin{pf}
The result is clear from the identities
\begin{eqnarray*}
\bbe\varphi(X_{\Delta})&=&\varphi(c)\bbp[
X_{\Delta}\geq c]+\int_{c}^{\infty}
\varphi'(u)\bbp[X_{\Delta}\geq u]\,\mathrm{d}u,\\
\int\varphi(x)\nu(\mathrm{d}x)&=&\varphi(c)\nu(
[c,\infty))+\int_{c}^{\infty}
\varphi'(u)\nu([u,\infty))\,\mathrm{d}u,
\end{eqnarray*}
which are standard consequences of Fubini's theorem.
\end{pf}

Our first result shows a central limit theorem for $\hat{s}(x)$
centered at $\bbe\hat{s}(x)$.
Let us remark that the fact that the Legendre polynomial $Q_{j}$ is not
constant for $j>0$ poses some difficulty since the relative position of
$x$ inside its class changes greatly with $m$.
\begin{lmma}\label{Lmma1CLT}
Under the notation and assumptions of Theorem \ref{TPWC}, it follows that
\[
\frac{c_{T}}{b_{k,m_{T}}(x)}\bigl(\hat{s}_{T}(x)-
\bbe\hat{s}_{T}(x)\bigr)
\,\ld\,\bar{\sigma} Z.
\]
\end{lmma}
\begin{pf}
We apply a central limit theorem version for row-wise independent
arrays of random variables (see, e.g.,
the corollary following \cite{Chung}, Theorem 7.1.2). Note that
\begin{eqnarray*}
S_{T}&:=&\frac{c_{T}}{b_{m_{T}}}\bigl(\hat{s}_{T}(x)-
\bbe\hat{s}_{T}(x)\bigr)
\\
&=&\frac{c_{T}}{Tb_{m_{T}}}
\sum_{i}
\sum_{j=0}^{k}\widetilde{\varphi}_{j,T}(x)
\{
\widetilde{\varphi}_{j,T}(X_{t^{i}_{T}}-X_{t^{i-1}_{T}})
-\bbe\widetilde{\varphi}_{j,T}(X_{\Delta_{T}^{i}})
\},
\end{eqnarray*}
where $\widetilde{\varphi}_{j,T}(\cdot)$ is of the form
\[
{\sqrt{\frac{2j+1}{b_{T}-a_{T}}}}
Q_{j}\biggl(\frac{2\cdot-(a_{T}+b_{T})}{b_{T}-a_{T}}
\biggr)\mathbf{1}_{[a_{T},b_{T})}(\cdot)
\]
with $a_{T},b_{T}$ such that
$x\in[a_{T},b_{T})$ and
$b_{T}-a_{T}=(b-a)/m_{T}$.
In that case,
$\bar{\sigma}_{T}^{2}:=\operatorname{Var} S_{T}$ is given by
%
\begin{equation}\label{Var2}
\bar{\sigma}_{T}^{2}:=
\frac{c^{2}_{T}}{T^{2}b_{m_{T}}^{2}}
\sum_{i}\sum_{j_{1},j_{2}=0}^{k}
\widetilde{\varphi}_{j_{1},T}(x)\widetilde{\varphi}_{j_{2},T}(x)
\operatorname{Cov}(\widetilde{\varphi}_{j_{1},T}(X_{\Delta_{T}^{i}}),
\widetilde{\varphi}_{j_{2},T}(X_{\Delta_{T}^{i}})),
\end{equation}
where we have used $\Delta^{i}_{T}:=t^{i}_{T}-t^{i-1}_{T}$.
Let us analyze the above covariances, scaled by $\Delta_{T}^{i}$.
First, applying
Lemma \ref{EstDstnc}, (\ref{EstimateLegendre}) and Proposition \ref
{UnfrmCvrgnc},
there exists a $t_{0}>0$ and $K>0$ such that
whenever $\Delta<t_{0}$,
\[
\biggl|\frac{1}{\Delta}
\bbe\widetilde{\varphi}_{j_{1},T}(X_{\Delta})
\widetilde{\varphi}_{j_{2},T}(X_{\Delta})
-\int\widetilde{\varphi}_{j_{1},T}(y)
\widetilde{\varphi}_{j_{2},T}(y)s(y)\,\mathrm{d}y\biggr|
\leq \frac{K\Delta}{b_{T}-a_{T}}.
\]
Similarly, using the additional fact that
$|\int\widetilde{\varphi}_{j,T}(y)s(y)\,\mathrm{d}y|
\leq \|s\|$, there exists a $t_{0}>0$ and $K>0$ such that
whenever $\Delta<t_{0}$,
\[
\biggl|\frac{1}{\Delta}
\bbe\widetilde{\varphi}_{j_{1},T}(X_{\Delta})
\bbe\widetilde{\varphi}_{j_{2},T}(X_{\Delta})
\biggr|\leq
K\Delta.
\]
Thus, using assumption (iii) of Theorem \ref{TPWC},
we have
\[
\frac{1}{\Delta_{T}^{i}}\operatorname{Cov}(\widetilde{\varphi
}_{j_{1},T}(X_{\Delta_{T}^{i}}),
\widetilde{\varphi}_{j_{2},T}(X_{\Delta_{T}^{i}}))
=\mathrm{o}_{T}(1)+\int\widetilde{\varphi}_{j_{1},T}(x)
\widetilde{\varphi}_{j_{2},T}(y)s(y)\,\mathrm{d}y,
\]
where $\mathrm{o}_{T}(1)\rightarrow 0$\vspace*{-2pt} uniformly in $i$ as $T\rightarrow
\infty$. Thus, in view of the fact that $b_{m_{T}}\geq 1$, (\ref
{EstimateLegendre}) and assumption (ii) of Theorem \ref{TPWC}, we have
$\bar{\sigma}_{T}^{2}-\hat{\sigma}_{T}^{2}
\stackrel{T\rightarrow\infty}{\longrightarrow}0$, where
\[
\hat{\sigma}_{T}^{2}:=
\frac{c^{2}_{T}}{T b_{m_{T}}^{2}}
\sum_{j_{1},j_{2}=0}^{k}
\widetilde{\varphi}_{j_{1},T}(x)\widetilde{\varphi}_{j_{2},T}(x)
\int\widetilde{\varphi}_{j_{1},T}(y)
\widetilde{\varphi}_{j_{2},T}(y)s(y)\,\mathrm{d}y.
\]
Next, the continuity of $s$ at $x$, assumption (ii) of Theorem \ref
{TPWC} and
the fact that the support of
$\widetilde{\varphi}_{j,T}$ contains $x$
and shrinks to $0$ collectively yield that
\[
\lim_{T\rightarrow\infty}
\frac{c^{2}_{T}}{Tb_{m_{T}}^{2}}
\sum_{j_{1},j_{2}=0}^{k}
\widetilde{\varphi}_{j_{1},T}(x)\widetilde{\varphi}_{j_{2},T}(x)
\int\widetilde{\varphi}_{j_{1},T}(y)
\widetilde{\varphi}_{j_{2},T}(y)\bigl(s(y)-s(x)\bigr)\,\mathrm{d}y=0.
\]
This implies that
$
\lim_{T\rightarrow\infty}
\hat{\sigma}_{T}^{2}=\lim_{T\rightarrow\infty}
\bar{\sigma}_{T}^{2}
=s(x)/(b-a),
$
in view of condition (ii) and the definition of
$b_{k}$.
Finally, we consider the ``standardized'' sum
$
Z_{T}:=S_{T}/\bar{\sigma}_{T}.
$
By the corollary following \cite{Chung}, Theorem 7.1.2,
$Z_{T}$ will converge to $\calN(0,1)$ because
\begin{eqnarray*}
&&\sup_{i}\frac{c_{T}}{T \bar{\sigma}_{T} b_{m_{T}}}
\sum_{j=0}^{k}|\widetilde{\varphi}_{j,T}(x)
\widetilde{\varphi}_{j,T}(X_{t^{i}_{T}}-X_{t^{i-1}_{T}})|\\
&&\quad\leq
\frac{c_{T}m_{T}}{T\bar{\sigma}_{T}
b_{m_{T}}(b-a)}\rightarrow 0
\end{eqnarray*}
as $T\rightarrow \infty$, in view of assumptions (i)--(ii) and the
fact that $ b_{m}\geq 1$. This implies the proposition since
$\bar{\sigma}_{T}^{2}\rightarrow s(x)(b-a)^{-1}$.
\end{pf}

The last step is to estimate the rate of convergence of the bias term.
\begin{lmma}\label{LRCB}
Under the notation and assumptions of Theorem \ref{TPWC},
$\bbe\hat{s}_{T}(x)-s(x)=\mathrm{o}(b_{m_{T}}/c_{T})$ as $T\to\infty$
for any fixed $x\in(a,b)$ such that $s(x)>0$.
\end{lmma}
\begin{pf}
We use the same notation as in the proof of Lemma \ref{Lmma1CLT}.
Obviously,
\[
\frac{c_{T}}{b_{m_{T}}}
|\bbe\hat{s}_{T}(x)-s(x)|
\leq\frac{1}{T}\sum_{i}\Delta_{T}^{i}
A_{T}(\Delta_{T}^{i}),
\]
where
\[
A_{T}(\Delta):=
\frac{c_{T}}{b_{m_{T}}}\Biggl|
\frac{1}{\Delta}
\sum_{j=0}^{k}\widetilde{\varphi}_{j,T}(x)
\bbe\widetilde{\varphi}_{j,T}(X_{\Delta})
-s(x)
\Biggr|.\vadjust{\goodbreak}
\]
It then suffices to show that
$\max_{i} A_{T}(\Delta^{i}_{T})\rightarrow 0$
as $T\rightarrow\infty$.
Note that
\begin{eqnarray*}
A_{T}(\Delta)&\leq&
\frac{c_{T}}{b_{m_{T}}}\Biggl|
\sum_{j=0}^{k}\widetilde{\varphi}_{j,T}(x)
\biggl\{
\frac{1}{\Delta}
\bbe\widetilde{\varphi}_{j,T}(X_{\Delta})
-\int\widetilde{\varphi}_{j,T}(y) s(y)\,\mathrm{d}y
\biggr\}
\Biggr|\\
 &&{}+
\frac{c_{T}}{b_{m_{T}}}\Biggl|
\int\sum_{j=0}^{k}\widetilde{\varphi}_{j,T}(x)
\widetilde{\varphi}_{j,T}(y) \bigl(s(y)-s(x)\bigr)\,\mathrm{d}y
\Biggr|,
\end{eqnarray*}
where we have used the fact that
$\int\widetilde{\varphi}_{j,T}(y) \,\mathrm{d}y=
\delta_{0}(j)$. We shall show that each of the two terms on the
right-hand side of the
above inequality, which we denote $A_{T}^{1}(\Delta)$ and
$A_{T}^{2}$, respectively, vanish as $T\to \infty$.
Using (\ref{EstimateLegendre}), Lemma \ref{EstDstnc} and Proposition
\ref{UnfrmCvrgnc},
there exist a $K>0$ and $T_{0}>0$ such that,
for $T>T_{0}$,
\begin{eqnarray*}
A_{T}^{1}(\Delta^{i}_{T})&\leq&
K \frac{c_{T}\Delta^{i}_{T}}{b_{m_{T}}(b_{T}-a_{T})}\\
&\leq& K \frac{c_{T}m_{T}\bar{\pi}_{T}}{(b-a)}\rightarrow 0
\end{eqnarray*}
as $T\rightarrow\infty$, due to (i)--(iii). To deal with the term
$A_{T}^{2}$,
we treat the two cases $\alpha=1$ and $\alpha>1$ separately. Suppose
that $\alpha=1$.
Using the Cauchy--Schwarz inequality twice (for summation and for the
integral) and the fact that
$\sum_{j=0}^{k} \widetilde{\varphi
}^{2}_{j,T}(x)=b^{2}_{m_{T}}(x)/(b_{T}-a_{T})$, we have
\[
A_{T}^{2}\leq\frac{c_{T}}{\sqrt{b_{T}-a_{T}}}
\Biggl\{\sum_{j=0}^{k}\int_{a_{T}}^{b_{T}} \bigl(s(y)-s(x)\bigr)^{2}
\,\mathrm{d}y\Biggr\}^{1/2}
\leq K c_{T}(b_{T}-a_{T})
\]
for some constant $K<\infty$. In light of assumption (iv) of Theorem
\ref{TPWC}, $A_{T}^{2}\stackrel{T\rightarrow\infty}{\longrightarrow}0$.
Let us now assume that $\alpha>1$. We first note that
\[
\int\sum_{j=0}^{k}\widetilde{\varphi}_{j,T}(x)
\widetilde{\varphi}_{j,T}(y) (y-x)^{j'}\,\mathrm{d}y=0
\]
for $j'=1,\ldots,k$. This is because the left-hand side is $p^{\bot
}(x)$, where
$p^{\bot}(y)$ is the orthogonal projection of the function $p(y):=(y-x)^{j'}$
on $\calS_{k,m_{T}}$ and, clearly,
$p^{\bot}(x)=p(x)=0$.
Also, by Taylor's theorem,
\[
s(y)-s(x)=\sum_{j'=1}^{r}
\frac{s^{(j')}(x)}{j'!}(y-x)^{j'}
+\int_{x}^{y} \bigl(s^{(r)}(v)-s^{(r)}(x)\bigr)
\frac{(y-v)^{r-1}}{(r-1)!}\,\mathrm{d}v,
\]
where $r:=\lfloor\alpha\rfloor$,
the largest integer that is (strictly) smaller than $\alpha$.
Since $k\geq \alpha-1$, we have that $k\geq r$ and
\begin{eqnarray*}
&&\int\sum_{j=0}^{k}\widetilde{\varphi}_{j,T}(x)
\widetilde{\varphi}_{j,T}(y) \bigl(s(y)-s(x)\bigr)\,\mathrm{d}y
\\
 &&\quad=\int\sum_{j=0}^{k}\widetilde{\varphi}_{j,T}(x)
\widetilde{\varphi}_{j,T}(y)\int_{x}^{y} \bigl(s^{(r)}(v)-s^{(r)}(x)\bigr)
\frac{(y-v)^{r-1}}{(r-1)!}\,\mathrm{d}v \,\mathrm{d}y.
\end{eqnarray*}
Again applying the Cauchy--Schwarz inequality twice (for summation and
for the integral),
we have
\begin{eqnarray*}
A_{T}^{2}
&\leq&
\frac{c_{T}}{b_{m_{T}}}
\sum_{j=0}^{k}|\widetilde{\varphi}_{j,T}(x)|
\biggl|\int
\widetilde{\varphi}_{j,T}(y)\int_{x}^{y} \bigl(s^{(r)}(v)-s^{(r)}(x)\bigr)
\frac{(y-v)^{r-1}}{(r-1)!}\,\mathrm{d}v \,\mathrm{d}y\biggr|\\
&\leq&
\frac{c_{T}}{\sqrt{b_{T}-a_{T}}}
\Biggl\{\sum_{j=0}^{k}
\int_{a_{T}}^{b_{T}}\biggl\{
\int_{x}^{y} \bigl(s^{(r)}(v)-s^{(r)}(x)\bigr)
\frac{(y-v)^{r-1}}{(r-1)!}\,\mathrm{d}v \biggr\}^{2}\,\mathrm{d}y\Biggr\}^{1/2}.
\end{eqnarray*}
Finally, by the H\"older condition (\ref{HolderCond}),
$
A_{T}^{2}\leq K c_{T}m_{T}^{-\alpha}\stackrel{T\rightarrow
\infty}{\longrightarrow}0.
$
\end{pf}

\section{Proofs of the uniform central limit theorem}\label{SectAppdxC}

In this section, we show the results of Section \ref{SectCBnds}. We
recall that the estimators $\hat{s}_{T}^{n}$ are based on
observation of the process at evenly-spaced times $\pi
_{T}^{n}\dvtx t_{0}=0<\cdots<t_{n}=T$. The time span between observations
is $\delta^{n}: =\delta^{n}_{T}:=T/n$.

Let us first remark that under the assumption that $\sigma\neq 0$ or
$\nu(\bbr)=\infty$,
the distribution $F_{t}(x)$ is continuous for all $t>0$ (see \cite
{Sato}, Theorem 27.4). In particular,
$\{F_{\delta^{n}}(X_{t_{i}}-X_{t_{i-1}})\}_{i\leq n}$ is necessarily a
random sample of uniform random variables and, hence, $Z_{n}^{0}$ of
(\ref{USEP}) is indeed the standardized empirical process of a uniform
random sample. Also, note that
\[
Z_{n}^{0}(F_{\delta^{n}}(x))=n^{1/2}\{F^{n}(x)-F_{\delta
^{n}}(x)\}\qquad \forall x\in\bbr,
\]
where $F^{n}:=F_{T}^{n}$ is the empirical process of $\{X_{t_{i}}-X_{t_{i-1}}\dvtx i=0,\ldots,n\}$.
The following transformation will be useful in the sequel:
\begin{eqnarray*}
\mathcal{L}(x;m,\kappa,H)
&=&\kappa\sum_{i=1}^{m}\sum_{j=0}^{k} \hat{\varphi}_{i,j}(x)
\biggl\{\hat{\varphi}_{i,j}(x_{i})\bigl(
H(x_{i})-H(x_{i-1})\bigr)\\
&&{}\hspace*{71pt} -\int_{x_{i-1}}^{x_{i}}\hat{\varphi}'_{i,j}(u)
\bigl(H(u)-H(x_{i-1})\bigr)\,\mathrm{d}u\biggr\},
\end{eqnarray*}
where $\hat{\varphi}_{i,j}$ is the basis element in (\ref
{BasisSplines}) and $H\dvtx \bbr\to\bbr$ is a locally integrable function.
Note that if $H$ is a function of bounded variation, then
\[
\mathcal{L}(x;m,\kappa,H) ={\kappa\sum_{i=1}^{m}\sum
_{j=0}^{k} \hat{\varphi}_{i,j}(x)\int_{x_{i-1}}^{x_{i}} \hat\varphi
_{i,j}(u)\,\mathrm{d} {H} (u)} .
\]
The following estimate follows easily from (\ref{EstimateLegendre}):
%
\begin{equation}\label{EqEstmTrnsf}
\sup_{x\in[a,b]}|\mathcal{L}(x;m,\kappa,H)|\leq
K\cdot\kappa\cdot m\cdot\omega\biggl(H;[a,b],\frac{b-a}{m}\biggr),
\end{equation}
where $K$ is a constant (depending only on $k$) and $\omega$ is the
modulus of continuity of $H$ defined by
\[
\omega(H;[a,b],\delta)=\sup\{|H(u)-H(v)|\dvtx u,v\in[a,b], |u-v|<\delta
\}.
\]

Let us write the estimator (\ref{ProjEstm}) in terms of
${F}^{n}_{T}$ as follows:
\begin{equation}\label{EqProjEstmSplines2}
\hat{s}^{n}_{T}(x):=\sum_{i=1}^{m}\sum_{j=0}^{k}
\hat\beta^{\pi_{T}^{n}}(\hat{\varphi}_{i,j})\hat{\varphi}_{i,j}(x)
=\mathcal{L}\biggl(x;m,{\frac{n}{T}},{F}^{n}_{T}(\cdot)\biggr).
\end{equation}
Note that $\bbe\hat{s}^{n}_{T}(x)$ admits a similar expression with
${F}^{n}_{T}$ replaced by
${F}_{\delta^{n}_{T}}$. Thus, it follows that a.s.
\begin{equation}\label{EqEmpExpr1}
Y^{n}_{T}(x)
:=\hat{s}^{n}_{T}(x)-\bbe\hat{s}^{n}_{T}(x)=
\calL(x;m,n^{1/2}T^{-1},Z_{n}^{0}(F_{\delta^{n}}(\cdot)))
\end{equation}
for all $x$. As was explained in Section \ref{SectCBnds}, one of the
key ideas of the approach of Bickel and Rosenblatt \cite{BicRos}
consists of approximating $Z_{n}^{0}$ by a Brownian bridge $Z^{0}$. To
this end, we use the following result, which follows from the Koml\'os, Major
and Tusn\'ady construction \cite{KMT}.
\begin{thrm}\label{ThrmApprxBrill}
There exists a probability space $(\widetilde\Omega,\widetilde\calF
,\widetilde\bbp)$, equipped with a standard Brownian motion $\widetilde
{Z}$, on which one can construct a version $\widetilde{Z}_{n}^{0}$ of
$Z_{n}^{0}$ such that
\[
{\|\widetilde{Z}_{n}^{0}-\widetilde{Z}^{0}\|_{[0,1]} =
\mathrm{O}_{p}(n^{-1/2}\log n)},
\]
where $\widetilde{Z}^{0}(x):= \widetilde{Z}(x)-x \widetilde{Z}(1)$ is
the corresponding Brownian bridge.
\end{thrm}

Since we are looking for the asymptotic distribution of $\sup
_{x}|Y_{T}^{n}(x)|$, properly scaled and centered, we can work with
the process $\widetilde{Z}_{n}^{0}$ instead of $Z_{n}^{0}$.
Thus, with some abuse of notation, we drop the tilde in all of the
processes of Theorem \ref{ThrmApprxBrill}. The following is an easy
estimate. Again abusing notation, the process $ _{0}Y_{T}^{n}$ in
the following lemma
is actually the process resulting from replacing $Z_{n}^{0}(F_{\delta
^{n}}(\cdot))$ in (\ref{EqEmpExpr1}) by
$\widetilde{Z}_{n}^{0}(F_{\delta^{n}}(\cdot))$.
\begin{lmma}\label{LEE3}
Let $_{0}Y_{T}^{n}(x)=
\calL(x;m,n^{1/2}T^{-1},Z^{0}(F_{\delta^{n}}(\cdot)))$. It
then follows that
$\| _{0}Y_{T}^{n}- Y_{T}^{n}\|_{[a,b]}=\mathrm{O}_{p}(m \log n/T
)$ as $n\to\infty$.
\end{lmma}
\begin{pf}
Clearly,
$
\omega(H;[a,b],\delta)\leq2 \| H\|_{[a,b]}
$
for any process $H$. Thus, we get the result from~(\ref{EqEstmTrnsf})
and Theorem \ref{ThrmApprxBrill}.
\end{pf}

As in \cite{BicRos},
our approach is to devise successive approximations of
$_{0}Y_{T}^{n}(x)$, denoted by $_{1}Y_{T}^{n},\ldots,
 _{N}Y_{T}^{n}$, such that the asymptotic distribution of the supremum
$
\sup_{x\in[a,b]}|  _{N}Y_{T}^{n}(x)|,
$
properly centered and scaled by certain constants $b_{T}^{n}$ and
$a_{T}^{n}$, is easy to determine and such that the error of the
successive approximations is negligible when multiplied by
$a^{n}_{T}$. We proceed to carry out this program.

First, note that since a Brownian bridge satisfies $\{Z^{0}(x)\}_{x\leq
 1}\,\ed\,\{Z^{0}(1-x)\}_{x\leq 1}$,
we have
\[
\{_{0}Y_{T}^{n}(x)\}_{x\in[a,b]}\,\ed\, \{ _{1}Y_{T}^{n}(x)\}
_{x\in[a,b]},
\]
where $_{1}Y_{T}^{n}(x):=
\calL(x;m,n^{1/2}T^{-1},Z^{0}(\bar{F}_{\delta^{n}}(\cdot))$ and
$\bar{F}:=1-F$. The following is our first estimate.
\begin{lmma}\label{LEE2}
Suppose that the assumptions of Proposition \ref{UnfrmCvrgnc} are
satisfied. There exist constants $K$ and $t_{0}>0$ such that if
$T/n<t_{0}$, then
\[
_{2}Y_{T}^{n}(x)=
\calL(x;m,n^{1/2}T^{-1},Z(\bar{F}_{\delta^{n}}(\cdot)))
\]
is such that
\[
\| _{1}Y_{T}^{n}-  {_{2}}Y_{T}^{n}\|_{[a,b]}\leq
K n^{-1/2}\biggl(\frac{mT}{n}\vee1\biggr)|Z(1)|
\]
for a constant $K<\infty$.
\end{lmma}
\begin{pf}
Clearly,
\[
_{2}Y_{T}^{n}(x)-  {_{1}}Y_{T}^{n}(x)=
\calL(x;n,T,m,n^{1/2}T^{-1},Z(1) \bar{F}_{\delta^{n}}(\cdot)).
\]
Thus, by (\ref{EqEstmTrnsf}),
\[
\| _{1}Y_{T}^{n}-  {_{2}}Y_{T}^{n}\|_{[a,b]}\leq
K \frac{m n^{1/2}}{T} \omega( \bar{F}_{\delta
^{n}};[a,b],d_{m}) |Z(1)|,
\]
where $d_{m}=(b-a)/m$. In view of Proposition \ref{UnfrmCvrgnc}, for
$n$ and $T$ such that
$T/n<t_{0}$, there are constants $k$ and $k'$ such that
\[
|\bar{F}_{\delta^{n}}(u)- \bar{F}_{\delta^{n}}(v)|\leq2k (\delta
^{n})^{2} +
2k' \delta^{n} m^{-1},
\]
provided that $u,v\in[a,b]$ and $|v-u|<d_{m}$.
\end{pf}

Let us now work with $_{2}Y_{T}^{n}$. Because of the self-similarity
of the Brownian motion, we have that
\[
\{_{2}Y_{T}^{n}(x)\}_{x\in[a,b]}\,\ed\, \{ _{3}Y_{T}^{n}(x)\}
_{x\in[a,b]},
\]
where
\[
_{3}Y_{T}^{n}(x):=
\calL\biggl(x;m,T^{-1/2},Z\biggl(\frac{1}{\delta^{n}}\bar{F}_{\delta^{n}}
(\cdot)\biggr)\biggr).
\]

The following estimate results from L\'evy's modulus of continuity theorem.
\begin{lmma}\label{LEE}
Let $_{4}Y_{T}^{n}(x)=
\calL(x;m,T^{-1/2},Z(\int_{\cdot}^{\infty} s(u)\,\mathrm{d}u
))$. If $T_{n}$ is such that
$\delta^{n}:=\frac{T_{n}}{n}\rightarrow 0$, then, for $n$ large enough,
\[
\| _{3}Y_{T_{n}}^{n}-  {_{4}}Y_{T_{n}}^{n}\|_{[a,b]}\leq
m\cdot \mathrm{O}_{p}\biggl({n}^{-1/2} \log^{1/2} \frac{n}{T_{n}}\biggr)
\]
for a constant $K<\infty$.
\end{lmma}
\begin{pf}
It is not hard to see that there exists a constant $K$ such that
\[
\| _{3}Y_{T}^{n}-  {_{4}}Y_{T}^{n}\|
\leq K T^{-1/2} m \sup_{x\in[a,b]}\biggl| Z\biggl(\frac{1}{\delta
^{n}}\bar{F}_{\delta^{n}}(x)\biggr) - Z\biggl(\int_{x}^{\infty}
s(u)\,\mathrm{d}u\biggr)\biggr|.
\]
By Proposition \ref{UnfrmCvrgnc}, there exist constants $k>0$ and
$t_{0}>0$ such that
for all $0<\delta<t_{0}$,
%
\begin{equation}
\sup_{y\in D} \biggl|
\frac{1}{\delta}\bbp[X_{\delta}\geq y]-
\nu([y,\infty))\biggr|
< k  \delta.
\end{equation}
Thus, there exists a constant $K>0$ such that, for large enough $n$,
\[
\| _{3}Y_{T}^{n}-  {_{4}}Y_{T}^{n}\|
\leq
K {n}^{-1/2} m \log^{1/2}\frac{n}{T_{n}}  \qquad  \mbox{a.s.}
\]
\upqed\end{pf}

We now note that
\[
\biggl\{Z\biggl(\int_{x}^{\infty}s(u) \,\mathrm{d}u\biggr)\biggr\}_{x\in[a,b]}
\,\ed\, \biggl\{\int_{x}^{\infty} s^{1/2} (u) \,\mathrm{d} Z(u)\biggr\}_{x\in[a,b]}
\]
and, hence,
\[
\{_{4}Y_{T}^{n}(x)\}_{x\in[a,b]}\,\ed\, \{ _{5}Y_{T}^{n}(x)\}
_{x\in[a,b]},
\]
where
\[
_{5}Y_{T}^{n}(x):=
\calL\biggl(x;m,T^{-1/2},\int_{\cdot}^{\infty}s^{1/2}(u) \,\mathrm{d}Z(u)\biggr).
\]
Using integration by parts, one can simplify $_{5}Y_{T}^{n}(x)$ as
follows:
\[
_{5}Y_{T}^{n}(x)=T^{-1/2}\sum_{i=0}^{m}\sum_{j=0}^{k} \hat{\varphi}_{i,j}(x)
\int_{x_{i-1}}^{x_{i}}s^{1/2}(u) \hat{\varphi}_{i,j}(u) \,\mathrm{d}Z(u).
\]

The following is the last estimate.
\begin{lmma}\label{ResErrorBnd}
Suppose that the Assumptions \ref{CndLevyDnsty} in Section \ref
{SectCBnds} hold true.
Let
\[
_{6}Y_{T}^{n}(x):=(b-a)^{1/2}T^{-1/2}
\sum_{i=0}^{m}\sum_{j=0}^{k} \hat{\varphi}_{i,j}(x)
\int_{x_{i-1}}^{x_{i}}\hat{\varphi}_{i,j}(u) \,\mathrm{d}Z(u).
\]
There then exists a random variable $M$ such that
\[
\| _{6}Y_{T}^{n}(\cdot)-
(b-a)^{1/2}s^{-1/2}(\cdot) _{5}Y_{T}^{n}(\cdot)\|
\leq M T^{-1/2}.
\]
\end{lmma}
\begin{pf}

Let $q(x)=s^{1/2}(x)$ and $c=(b-a)^{1/2}$. Using integration by parts,
we have
\begin{eqnarray*}
H_{i,j}(x)&:=&s^{-1/2}(x)\int_{x_{i-1}}^{x_{i}}s^{1/2}(u) \hat{\varphi
}_{i,j}(u) \,\mathrm{d}Z(u)-\int_{x_{i-1}}^{x_{i}} \hat{\varphi}_{i,j}(u) \,\mathrm{d}Z(u)\\
 &=&q^{-1}(x)\bigl\{\hat{\varphi}_{i,j}(x_{i})
\bigl(q(x_{i})-q(x)\bigr) Z(x_{i}) -
\hat{\varphi}_{i,j}(x_{i-1})\bigl(q(x_{i-1})-q(x)
\bigr)Z(x_{i-1})\bigr\}\\
&&{}-q^{-1}(x)\int_{x_{i-1}}^{x_{i}}\bigl\{\hat{\varphi
}_{i,j}'(u)\bigl(
q(u)-q(x)\bigr)-\hat{\varphi}_{i,j}(u)q'(u)\bigr\}Z(u)\,\mathrm{d}u.
\end{eqnarray*}
Since $q^{-1}(\cdot)$ and $q'(\cdot)$ are bounded on $[a,b]$, there
exists a constant $K$ such that
\[
\sup_{x\in[x_{i-1},x_{i}]}|H_{i,j}(x)|\leq K
m^{-1/2} \sup_{u\in[x_{i-1},x_{i}]}|Z(u)|.
\]
Thus,
\begin{eqnarray*}
\| _{6}Y_{T}^{n}(\cdot)-
cs^{-1/2}(\cdot) _{5}Y_{T}^{n}(\cdot)\|&
\leq&\biggl(\frac{T}{b-a}\biggr)^{-1/2} \sum_{i=0}^{m}\sum_{j=0}^{k}\sup
_{x\in[x_{i-1},x_{i}]}|H_{i,j}(x)\hat{\varphi}_{i,j}(x)|\\
&\leq&
K
T^{-1/2} \sup_{u\in[a,b]}|Z(u)|.
\end{eqnarray*}
\upqed\end{pf}

{The latter approximation, $_{6}Y_{T}^{n}$, is simple enough to try
determining its asymptotic distribution (appropriately centered and
scaled). Indeed, }
%
\begin{equation}\label{EqDstrbMxm}
M(T,n,m):=\sup_{x\in[a,b]} |_{6}Y_{T}^{n}(x)|\,\ed\,
T^{-1/2}m^{1/2}\max_{1\leq j\leq m} \bigl\{ \zeta_{m}^{(k)}\bigr\},
\end{equation}
where $\{\zeta_{j}^{(k)}\}_{i}$ are independent copies of the
r.v. $\zeta^{(k)}$ defined in (\ref{EqKeySprm}).\vspace*{-2pt} The following result
obtains the asymptotic distributions of
$
\bar{M}_{m}:=\max_{1\leq j\leq m}\{\zeta_{j}^{(k)}\}
$
for the cases $k=0$ and $k=1$.
\begin{lmma}
Let $a_{n}$ and $b_{n}$ be as in \textup{(\ref{NrmlzCnst1})--(\ref
{NrmlzCnst2})}. The following limits then hold:
\begin{eqnarray}\label{EqLimitMaxima0}
\lim_{m\rightarrow\infty}\bbp\biggl(
\max_{1\leq j\leq m} \bigl\{ \zeta_{m}^{(0)}\bigr\}\leq
\frac{{y}}{a_{m_{n}}}+b_{m_{n}}\biggr)&=&\mathrm{e}^{-2\mathrm{e}^{-{y}}},\\
\label{EqLimitMaxima20}
\lim_{m\rightarrow\infty}\bbp\biggl( 2^{-1} \max_{1\leq j\leq m}\bigl \{
\zeta_{m}^{(1)}\bigr\}\leq\frac{{y}}{a_{m_{n}}}+b_{m_{n}}
\biggr)&=&\mathrm{e}^{-4\mathrm{e}^{-{y}}}
\end{eqnarray}
for all ${y}\in\bbr_{+}$.
\end{lmma}
\begin{pf}
The limit (\ref{EqLimitMaxima0}) follows from the well-known identity
%
\begin{equation}\label{BscLmt}
\lim_{m\rightarrow\infty}m\bigl(1-\Phi(u_{m}({y}))\bigr)=\mathrm{e}^{-{y}},
\end{equation}
where $\Phi$ is the normal distribution and $u_{m}(y)={y}/a_{m}+b_{m}$.
Indeed, for large enough $m$, the probability in (\ref{EqLimitMaxima0})
can be written as follows:
\[
\bigl(2\Phi(u_{m}({y}))-1\bigr)^{m}=
\biggl(1-\frac{2m(1-\Phi(u_{m}({y})))}{m}
\biggr)^{m}\longrightarrow \mathrm{e}^{-2\mathrm{e}^{-{y}}}.
\]
To handle the case $k=1$, we embed the problem into the theory of
multivariate extreme values (see, e.g., \cite{Galambos}). Consider
independent copies $\{\mathbf{V_{i}}\}_{i}$ of the following vector of
jointly standard Gaussian variables:
%
\begin{equation}\label{EqMultVect}
\mathbf{V}:=\biggl(\frac{1}{2}Z_{0}+\frac{\sqrt{3}}{2}Z_{1},\frac
{1}{2}Z_{0}-\frac{\sqrt{3}}{2} Z_{1}\biggr)'.
\end{equation}
Since $\zeta^{(1)}=|Z_{0}|+\sqrt{3}|Z_{1}|$, we can see that
\begin{eqnarray*}
&&\biggl\{{2^{-1}} \max_{1\leq j\leq m} \bigl\{ \zeta_{m}^{(1)}\bigr\}
\leq{\frac{{y}}{{a}_{m}}+{b}_{m}}\biggr\} \\
&&\quad=\Bigl\{\max_{i\leq m}\mathbf{V}_{i}\leq
{\bf\hat{a}}_{m}^{-1}{\mathbf{y}} + {\bf\hat{b}}_{m},
\min_{i\leq m}\mathbf{V}_{i}\geq
-{\bf\hat{a}}_{m}^{-1}{\mathbf{y}} - {\bf\hat{b}}_{m}\Bigr\},
\end{eqnarray*}
where {${\mathbf{y}}:=({y},{y})'$, ${\bf\hat{b}}_{m}:=({b}_{m},{b}_{m})'$,
${\bf\hat{a}}_{m}:=({a}_{m},{a}_{m})'$} and all operations are pointwise.
Then, (\ref{EqLimitMaxima20}) will follow from the following identity:
\begin{eqnarray}\label{EqMultLimit}
&&\lim_{m\rightarrow\infty}\bbp\Bigl(\max_{1\leq i\leq m}\mathbf{V}_{i}\leq
{\bf\hat{a}}_{m}^{-1}{\mathbf{y}} + {\bf\hat{b}}_{m},
\min_{1\leq i\leq m}\mathbf{V}_{i}\geq
-{\bf\hat{a}}_{m}^{-1}\mathbf{z} - {\bf\hat{b}}_{m}\Bigr)
\nonumber
\\[-8pt]
\\[-8pt]
\nonumber
&&\quad=
\mathrm{e}^{-\mathrm{e}^{-y_{1}}-\mathrm{e}^{-y_{2}}-\mathrm{e}^{-z_{1}}-\mathrm{e}^{-z_{1}}}\nonumber
\end{eqnarray}
for any ${\mathbf{y}}=({y}_{1},{y}_{2})'$ and ${\mathbf{z}}=(z_{1},z_{2})'$.
To show (\ref{EqMultLimit}), first note that the probability therein
can be written as
\[
A_{n}:=\bigl\{\bbp\bigl(-u_{n}(z_{1})\leq V_{1} \leq u_{n}(y_{1}),
-u_{n}(z_{2})\leq V_{2} \leq u_{n}(y_{2})\bigr)\bigr\}^{n},
\]
where $\mathbf{V}:=(V_{1},V_{2})'$ is defined in (\ref{EqMultVect}) and
$u_{n}(x):=x/a_{n}+b_{n}$.
Let
\[
{\bar\mathbf{F}_{n}(y,z;X,Y):=\bbp\bigl(X\geq u_{n}(y),Y\geq
u_{n}(z)\bigr),\qquad
\bar{F}_{n}(y;X):=\bbp\bigl(X\geq u_{n}(y)\bigr),}
\]
where $X$ and $Y$ represent random variables. We {recall} the
following results valid for any jointly normal variables $X$ and $Y$
and arbitrary $y$ and $z$ (see \cite{Galambos}, Example 5.3.1):
\[
{\lim_{n\rightarrow\infty}
n\bar\mathbf{F}_{n}(y,z;X,Y)=0,\qquad
\lim_{n\rightarrow\infty}
n\bar{ F}_{n}(y;X)=\mathrm{e}^{-y}.}
\]
Then, (\ref{EqMultLimit}) follows once we note that $A_{n}^{1/n}$ can
be written as follows:
\begin{eqnarray*}
A_{n}^{1/n}&=&1-\frac{1}{n}\{
n\bar{ F}_{n}(z_{1};V_{1})+n\bar{ F}_{n}(z_{2};V_{2})+n\bar{
F}_{n}(y_{1};-V_{1})+n\bar{ F}_{n}(y_{2};-V_{2})\\
&&\hspace*{26pt}{} -n\bar\mathbf{F}_{n}(z_{1},z_{2};V_{1},V_{2})
-n\bar\mathbf{F}_{n}(y_{1},z_{2};-V_{1},V_{2})
-n\bar\mathbf{F}_{n}(z_{1},y_{2};V_{1},-V_{2})\}.
\end{eqnarray*}
\upqed\end{pf}

In view of (\ref{EqDstrbMxm}), the following are easy consequences of
the above lemma:
\begin{eqnarray}\label{EqLimitMaxima}
\lim_{n\rightarrow\infty}\bbp\biggl(
T_{n}^{1/2} m_{n}^{-1/2}
\sup_{x\in[a,b]} |_{6}Y_{T_{n}}^{n}(x)|\leq
\frac{{y}}{a_{m_{n}}}+b_{m_{n}}\biggr)&=&\mathrm{e}^{-2\mathrm{e}^{-{y}}},\\
\label{EqLimitMaxima2}
\lim_{n\rightarrow\infty}\bbp\biggl(
{2^{-1}}T_{n}^{1/2} m_{n}^{-1/2}
\sup_{x\in[a,b]}| _{6}Y_{T_{n}}^{n}(x)|\leq\frac
{{y}}{a_{m_{n}}}+b_{m_{n}}\biggr)&=&\mathrm{e}^{-4\mathrm{e}^{-{y}}},
\end{eqnarray}
valid for all ${y}\in\bbr_{+}$, $T_{n}>0$ and $m_{n}$ such that
$m_{n}\rightarrow\infty$.
We are now ready to prove the main theorem of Section \ref{SectCBnds}:
\begin{pf*}{Proof of Theorem \protect\ref{TCB}}
The idea is to use the following simple observations. Let $\calL_{n}$
be a functional on $D[a,b]$ such that
%
\begin{equation}\label{EqLpschtzCond}
|\calL_{n}(\omega_{1})-\calL_{n}(\omega_{2})|\leq
M_{n}\|\omega_{1}-\omega_{2}\|
\end{equation}
and let $A_{n},B_{n}$ be processes with values on $D[a,b]$ such that
$\|A_{n}-B_{n}\|=\mathrm{o}_{p}(1/M_{n})$. Then, if $\calL_{n}(A_{n})$
converges in distribution to
$F$, $\calL_{n}(B_{n})$ will also converge to $F$. Throughout this
proof,\vspace{-5pt}
\[
\calL_{n}(\omega):=
a_{m_{n}}\biggl\{
\kappa  \cdot\frac{c}{d}\cdot\frac{T_{n}^{1/2}}{m_{n}^{1/2}}\cdot
\sup_{x\in[a,b]} | s^{-1/2}(x)\omega(x)|-b_{m_{n}}\biggr\},
\]
which satisfies the Lipschitz condition (\ref{EqLpschtzCond}) with
$M_{n}=\frac{\kappa c}{d} a_{m_{n}} T_{n}^{1/2}/m_{n}^{1/2}$.
From Lemma \ref{ResErrorBnd}, in order for
(\ref{EqLimitMaxima2}) to hold with $_{6}Y_{T_{n}}^{n}$ replaced by
$ _{5}Y_{T_{n}}^{n}$,
it suffices that
\[
\lim_{n\rightarrow\infty}\frac{T_{n}^{1/2}}{m_{n}^{1/2}}a_{m_{n}}
T_{n}^{-1/2} =
\lim_{n\rightarrow\infty}\biggl(\frac{2\log m_{n}}{m_{n}}\biggr)^{1/2}=0,
\]
which is obvious since $m_{n}\to\infty$. Since $ _{4}Y_{T_{n}}^{n}$
has the same law as $ _{5}Y_{T_{n}}^{n}$,
(\ref{EqLimitMaxima2}) also holds for $ _{4}Y_{T_{n}}^{n}$. In the
light of Lemma~\ref{LEE},
(\ref{EqLimitMaxima2}) will hold for $ _{3}Y_{T_{n}}^{n}$ (and,
hence, for
$ _{2}Y_{T_{n}}^{n}$ as well) since
\[
\lim_{n\rightarrow\infty}
\frac{T_{n}^{1/2}}{m_{n}^{1/2}}a_{m_{n}}
m_{n}{n}^{-1/2} \log^{1/2} \frac{n}{T_{n}}=
c \lim_{n\rightarrow\infty}
\biggl( m_{n}\log m_{n} \cdot\frac{T_{n}}{n} \log\frac{n}{T_{n}}
\biggr)^{1/2} =0,
\]
which follows from condition (ii) in the statement of Theorem \ref{TCB}.
Similarly, in view of Lemma~\ref{LEE2}, (\ref{EqLimitMaxima2}) will
hold for
$ _{1}Y_{T_{n}}^{n}$ (and hence, for
$ _{0}Y_{T_{n}}^{n}$ as well) since
\[
\lim_{n\rightarrow\infty}
\frac{T_{n}^{1/2}a_{m_{n}}
n^{-1/2}}{m_{n}^{1/2}}\biggl(\frac{m_{n}T_{n}}{n}\vee1\biggr) =0.
\]
Indeed,the above expression is upper bounded by
$
(\frac{T_{n}m_{n}}{n})^{1/2}
\frac{\log^{1/2} m_{n}}{m_{n}},
$
which converges to $0$ because of assumption (i) and the fact that
$m_{n}\to \infty$.
Finally, in the light of Lemma \ref{LEE3},
in order for (\ref{EqLimitMaxima2}) to hold for $Y_{T_{n}}^{n}$,
it suffices that
\[
\lim_{n\rightarrow\infty}
\frac{T_{n}^{1/2}}{m_{n}^{1/2}}a_{m_{n}}
\frac{m_{n}}{T_{n}} \log n=0,
\]
which follows from assumption (ii) in the statement of Theorem \ref{TCB}.
\end{pf*}

\begin{pf*}{Proof of Corollary \protect\ref{KRCB}}
Using the same reasoning as in the proof of Theorem \ref{TPWC}, it
turns out that
\[
\sup_{x\in[a,b]}|\bbe\hat{s}^{n}_{T_{n}}(x)-s(x)|
\leq K\biggl( \frac{m_{n}T_{n}}{n} \vee m_{n}^{-\alpha}\biggr)
\]
for an absolute constant $K$.
As in the proof of Theorem \ref{TCB}, to show (\ref{EqLimitMaxima3}),
it suffices that
\[
\lim_{n\rightarrow\infty} \frac{T_{n}^{1/2}}{m_{n}^{1/2}}a_{m_{n}}
\biggl( \frac{m_{n}T_{n}}{n} \vee m_{n}^{-\alpha}\biggr)=0,
\]
which holds in light of assumption (iii) in the statement of Corollary
\ref{KRCB}.
\end{pf*}
\end{appendix}

\section*{Acknowledgements}

The author's research was partially supported by NSF Grant No. DMS
0906919. The author is indebted to the referee and Editor for their
many suggestions that improved the paper considerably.
It is also a great pleasure to thank Professor David Mason for pointing
out the KMT inequality and for other important remarks. The author
would also like to thank Professor Jayanta Ghosh and participants of
the Workshop on Infinitely Divisible Processes (CIMAT A.C. March 2009)
for their helpful feedback.

\printhistory

\end{document}